\documentclass[11pt,leqno]{amsart}
\usepackage{amsmath,amsfonts,amssymb,amsthm,upref,bbm}
\usepackage{amscd}
\usepackage{enumitem}
\usepackage{needspace}
\usepackage{marginnote}

\usepackage[hidelinks]{hyperref}
\setlist[enumerate,1]{label=(\arabic*),leftmargin=2em,itemsep=0pt,parsep=0pt,topsep=3pt,partopsep=0pt}
\newtheoremstyle{compact}{5pt plus 1pt minus 1pt}{5pt plus 1pt minus 1pt}{\itshape}{}{\bfseries}{.}{.5em}{}
\theoremstyle{compact}
\AddToHook{cmd/normalsize/after}{%
    \setlength{\abovedisplayskip}{8pt plus 2pt minus 2pt}%
    \setlength{\belowdisplayskip}{8pt plus 2pt minus 2pt}%
    \setlength{\abovedisplayshortskip}{3pt plus 1pt}%
    \setlength{\belowdisplayshortskip}{5pt plus 1pt minus 1pt}%
}
\hypersetup{pdftitle={draft}}

\newtheorem{maintheorem}{Theorem}

\numberwithin{equation}{section}
\newtheorem{theorem}{Theorem}[section]
\newtheorem{corollary}[theorem]{Corollary}

\newtheorem{lemma}[theorem]{Lemma}
\newtheorem{proposition}[theorem]{Proposition}

\newtheorem{definition}[theorem]{Definition}
\newtheoremstyle{compactremark}{5pt plus 1pt minus 1pt}{5pt plus 1pt minus 1pt}{\normalfont}{}{\bfseries}{.}{.5em}{}
\theoremstyle{compactremark}
\newtheorem{example}[theorem]{Example}
\newtheorem{remark}[theorem]{Remark}
\theoremstyle{compact}

\DeclareMathOperator{\spn}{span}
\def\cspan{\overline{\operatorname{span}}}
\newcommand{\1}{\mathbbm{1}}
\newcommand{\rmd}{\mathrm{d}}
\newcommand{\N}{\mathbb{N}}
\newcommand{\R}{\mathbb{R}}
\newcommand{\cB}{\mathcal{B}}

\newcommand{\norm}[1]{\left\|#1\right\|}
\newcommand{\abs}[1]{\left|#1\right|}
\newcommand{\inp}[2]{\langle #1,#2\rangle}

\allowdisplaybreaks

\title[The Geometry of Banach space OVMs with cyclic dilations]{The Geometry of Banach space operator-valued measures with cyclic dilations}
\date{November 22, 2026}

\author{Daniel Freeman}
\address{Department of Mathematics and Statistics, St Louis University, St Louis, MO 63103, USA}
\email{daniel.freeman@slu.edu}

\author{Qianfeng Hu}
\address{School of Science, Hebei University of Technology, Tianjin, China}
\email{qianfenghu@hebut.edu.cn}

\author{Rui Liu}
\address{School of Mathematical Sciences and LPMC, Nankai University, Tianjin, 300071, PR China.}
\email{ruiliu@nankai.edu.cn}

\author{Jie Shen}
\address{School of Mathematical Sciences and LPMC, Nankai University, Tianjin, 300071, PR China.}
\email{1710064@mail.nankai.edu.cn}


\begin{document}

\begin{abstract}
In this paper, we investigate the geometric properties of operator-valued measures with cyclic dilations. Under a localization condition on the normal part of the bidual, we establish James-type characterizations of shrinking, bounded completeness, and reflexivity for projection-valued measures. 
These characterizations apply to spaces admitting no continuous Schauder frame, including $L_1[0,1]$. Using vector measure representations and interpolation, we prove that an operator-valued measure admitting a cyclic dilation has a reflexive cyclic dilation if and only if the underlying space is reflexive. Under the same cyclicity assumption, we also characterize the existence of cyclic dilations on spaces containing no $\ell_1$ and on KB-spaces. The dilation spaces can be chosen to be order-continuous K\"othe function spaces.
Without assuming cyclicity, we give a weak compactness criterion and local variation conditions for reflexive dilations; an example on $\ell_2$ shows that reflexivity alone is insufficient. For a separable space containing no $c_0$, a normalized operator-valued measure is induced by a continuous Schauder frame if and only if it admits a cyclic dilation and has $\sigma$-finite operator-norm variation. Finally, for measure quotients whose synthesis maps admit bounded linear right inverses, we construct compatible cyclic $L_p$-dilations and represent the quotient spaces as inductive limits, with applications to ReLU Barron spaces.

\end{abstract}

\maketitle

\section{Introduction}
Dilation theory is a paradigm by which we can exhibit operators as compressions of operators that are in some sense well-behaved, as well as view frames as compressions of bases~\cite{Pa2002,HanLarson2000}. Classical dilation theorems, including Sz.-Nagy's theorem on unitary dilations for contractions, Naimark's dilation theorem on orthogonal projection-valued dilations for positive regular operator-valued measures, and  $*$-homomorphism dilations due to Steinspring for completely positive maps on $C^*$-algebras, play an essential role in operator algebras. In these theorems, the dilation spaces are Hilbert spaces, and such dilations are called \textit{Hilbertian dilations}. However, an operator-valued measure acting on a Hilbert space need not have a Hilbertian dilation. From the dilation viewpoint, a Schauder frame can be viewed as a compression of a basis, which relates to bounded approximation properties of Banach spaces \cite{CHL1999} and leads to the study of dilation theory in the general setting of Banach spaces.

Han et al.~\cite{HanLarsonLiuLiu2014} established a Banach space dilation theory in which every operator-valued measure has a projection-valued dilation and every bounded linear map from a unital Banach algebra into $\mathcal B(X)$ admits a bounded homomorphism dilation. It is natural to ask whether the structure of the underlying spaces or maps can be preserved under dilations. Han et al.~\cite{HanLarsonLiu2018} proved that every operator-valued measure with bounded $p$-variation admits a projection-valued dilation with bounded $p$-variation. Furthermore, when the domain space is a von Neumann algebra, we often deal with normal (i.e., ultraweakly or $w^{*}$-continuous) linear maps. Thus it is natural to ask whether there exists a dilation that also preserves such continuity or normality of the linear maps. Han et al. obtained normal homomorphism dilations for $\ell^\infty$-spaces over purely atomic measure spaces \cite{HanLarsonLiuLiu2014} and $L^\infty$-spaces over finite measure spaces \cite{HanLarsonLiu2018}, and reduced the problem to countably additive projection-valued dilations of operator-valued quantum measures on projection lattices to acquire normal Jordan homomorphism dilations in the noncommutative setting \cite{HanHuLarsonLiu2024}. The geometric properties of the dilation spaces are a key ingredient in these constructions, and Banach space techniques are heavily involved. We therefore study cyclic dilations that preserve reflexivity or the absence of copies of $\ell_1$ or $c_0$. And we also obtain the Radon-Nikodým property under additional assumptions and use it to establish continuous frame representations.

Let $X$ be a nonzero real Banach space and $\mu\colon\Sigma\to\mathcal B(X)$ be a normalized, strongly countably additive operator-valued measure on $(\Omega,\Sigma)$ and $\mu(\Omega)=I_X$. A dilation system consists of a Banach space $Z$, bounded operators $T\colon X\to Z$ and $S\colon Z\to X$, and a normalized, strongly countably additive projection-valued measure $P$ on $Z$ such~that
\[
\mu(E)=SP(E)T\quad(E\in\Sigma),\qquad ST=I_X.
\]
Thus $P(E)P(F)=P(E\cap F)$, and $X$ is isomorphic to a complemented subspace of $Z$. Consequently, reflexivity of $X$ is necessary for a reflexive dilation. A projection-valued measure  $P$ on $Z$ is called \textit{cyclic} if there exists $u\in Z$ such that
\[
Z=\overline{\operatorname{span}}\{P(E)u:E\in\Sigma\}.
\]
Such a vector $u$ is called a cyclic vector for $P$. A dilation system $(Z,P,S,T)$ is called cyclic if $P$ is cyclic.

Cyclic projection-valued measures occur in familiar Banach spaces. If $(u_n)$ is an unconditional basis of $Z$ and $P$ is its coordinate projection-valued measure, then every $u=\sum_n a_nu_n\in Z$ with all $a_n\ne0$ is cyclic, since $P(\{n\})u=a_nu_n$. Likewise, on an order-continuous Banach function space $Z$, every $u\in Z$ that is nonzero almost everywhere is cyclic for the multiplication projection-valued measure $P(E)f=\mathbf{1}_Ef$. Cyclicity also arises naturally in the study of group frame generators~\cite{HanLarson2000,HanHuLiuLiu2026}. A projective unitary representation of a countable group $G$ with dense Bessel vectors is cyclic if and only if it admits a Parseval frame vector. Such a representation is unitarily equivalent to a subrepresentation of the corresponding left regular projective representation~\cite[Prop.~1.4 and Thm.~1.14]{HanHuLiuLiu2026}. 
In general, the minimal dilation space need not be reflexive even the underlying space is a Hilbert space, also see Remark~\ref{rem:scalar-minimal}. Here we consider cyclic dilations in Banach spaces and determine when the dilation space can be chosen reflexive, or can retain the absence of copies of $\ell_1$ or $c_0$ in the underlying space.


The theorem of James characterizes reflexivity of a Banach space with a Schauder basis in terms of shrinking and bounded completeness~\cite{AlbiacKalton2016}.
For unconditional bases, these properties are also characterized by the absence of copies of $\ell_1$ and $c_0$, respectively. For discrete Schauder frames, Beanland et al.~\cite{BeanlandFreemanLiu2015} proved that a frame has a shrinking associated basis if and only if it is shrinking, and has a reflexive associated space if and only if it is shrinking and boundedly complete. Beanland and Freeman~\cite{BeanlandFreeman2024} gave an explicit construction of shrinking associated bases for shrinking Schauder frames. Carando et al.~\cite{CarandoLassalleSchmidberg2011} obtained the corresponding characterizations for unconditional Schauder frames and atomic decompositions. Eisner and Freeman~\cite{EisnerFreeman2020} developed James-type results for continuous Schauder frames. We extend this correspondence to projection-valued measures under a localization condition on the normal part of the bidual.

For a projection-valued measure $P$ on $X$, its normal part $\mathcal N_P$ consists of those $\xi\in X^{**}$ for which $E\mapsto\xi(P(E)^*f)$ is countably additive for every $f\in X^*$. The localization condition~\eqref{eq:localization} requires that, for every $\xi\in\mathcal N_P$, there are measurable sets $E_n\uparrow\Omega$ such that $P(E_n)^{**}\xi\in J_XX$ for all $n$. Here $J_X\colon X\to X^{**}$ is the canonical embedding. The notions of shrinking and bounded completeness for $P$ are given in Section~3.

\Needspace{12\baselineskip}
\begin{maintheorem}[Theorem~\ref{thm:representation-james}]\label{intro:james}
Let $P$ be a normalized, strongly countably additive projection-valued measure on $X$ satisfying the localization condition~\eqref{eq:localization}. Then:
\begin{enumerate}
\item $P$ is shrinking if and only if $X$ contains no copy of $\ell_1$, and this is equivalent to $\mathcal N_P=X^{**}$.
\item $P$ is boundedly complete if and only if $X$ contains no copy of $c_0$. These conditions are equivalent to weak sequential completeness of $X$ and to $\mathcal N_P=J_XX$.
\item $X$ is reflexive if and only if $P$ is shrinking and boundedly complete.
\end{enumerate}
In particular, these assertions apply to cyclic and finitely generated representations.
\end{maintheorem}

This formulation includes spaces which admit neither an unconditional Schauder frame nor a continuous Schauder frame. For example, on $L_1[0,1]$, the multiplication measure $P(E)h=\1_Eh$ is cyclic and boundedly complete but not shrinking, although $L_1[0,1]$ admits no such frame (Example~\ref{ex:l1-james}). Thus the projection-valued formulation brings $L_1[0,1]$ within this James-type correspondence beyond the unconditional basis and continuous frame settings. The Orlicz heart $(\exp L_2)_0$ gives the opposite case: its multiplication measure is shrinking but not boundedly complete, and this space also admits no such frame (Example~\ref{ex:orlicz-james}).

The normal part identifies both properties within the inclusions $J_XX\subseteq\mathcal N_P\subseteq X^{**}$. Moreover, the localization condition has a three-space property for closed invariant subspaces (Proposition~\ref{prop:normal-coverage}). The characterizations for finitely generated Banach $C(K)$-modules and modules of finite multiplicity were obtained by Kitover and Orhon~\cite{KitoverOrhon2014,KitoverOrhon2017}. Our condition also applies to representations beyond finite multiplicity, including examples on spaces not isomorphic to Banach lattices (Example~\ref{ex:vector-sequences}).


We next use these characterizations to preserve geometric properties under cyclic dilations. Recall that a KB-space is a Banach lattice in which every increasing norm-bounded positive sequence converges in norm.

\Needspace{12\baselineskip}
\begin{maintheorem}[Theorem~\ref{thm:dilation-geometry}]\label{intro:dilation}
Let $\mu$ be a normalized, strongly countably additive operator-valued measure on $X$ admitting a cyclic dilation. In each assertion below, the dilation space $Z$ can be chosen to be an order-continuous K\"othe function space.
\begin{enumerate}
\item $X$ contains no copy of $\ell_1$ if and only if $\mu$ has a cyclic dilation $(Z,P,S,T)$ with $Z$ containing no copy of $\ell_1$ and $P$ shrinking.
\item $X$ contains no copy of $c_0$ if and only if $\mu$ has a cyclic dilation $(Z,P,S,T)$ with $Z$ a KB-space and $P$ boundedly complete.
\item $X$ is reflexive if and only if $\mu$ has a reflexive cyclic dilation.
\end{enumerate}
\end{maintheorem}

The dilation in each case represents the given measure $\mu$, with $\mu(E)=SP(E)T$ for every $E\in\Sigma$. To construct these dilation systems, we first show that $\mu$ admits a cyclic dilation if and only if its minimal $\alpha$-dilation is cyclic. In this case, the minimal dilation space is isomorphic to $L_1(m)$ for a countably additive $X$-valued measure $m$. With this representation, we apply the Davis--Figiel--Johnson--Pe\l czy\'nski construction~\cite{DFJP1974} to sets invariant under bounded measurable multipliers. Multiplier invariance and order continuity ensure that the multiplication projections remain strongly countably additive on the resulting spaces.

Without assuming cyclicity, we characterize the existence of a reflexive dilation by relative weak compactness of a multiplier orbit in the minimal dilation space (Theorem~\ref{prop:multiplier-orbit-criterion}). We also prove that reflexivity of $X$ is sufficient when every vector measure $E\mapsto\mu(E)x$ has finite variation on each set of a common countable measurable partition (Theorem~\ref{thm:variation-reflexive-dilation}). On the other hand, a normalized operator-valued measure on $\ell_2$ with bounded $2$-variation need not have a reflexive dilation (Example~\ref{ex:no-reflexive-dilation}). This example shows why an additional condition on the measure is needed. And we show that a normalized operator-valued measure on a separable Banach space containing no copy of $c_0$ is induced by a continuous Schauder frame on the same measurable space with a $\sigma$-finite representation measure if and only if it admits a cyclic dilation and has $\sigma$-finite operator-norm variation (Theorem~\ref{thm:continuous-frame-representation}).

Li, Li, and Han~\cite{LiLiHan2016} studied continuous framings for Banach spaces through their connection with $g$-K\"othe function spaces. Continuous Schauder frames give cyclic dilations, but a cyclic dilation need not arise from a continuous Schauder frame. The following result characterizes when a given measure has such a representation on its original measurable space.

\Needspace{12\baselineskip}
\begin{maintheorem}[Theorem~\ref{thm:continuous-frame-representation}]\label{intro:frame}
Let $X$ be a separable Banach space containing no copy of $c_0$, and let $\mu\colon\Sigma\to\mathcal B(X)$ be a normalized, strongly countably additive operator-valued measure. The following are equivalent:
\begin{enumerate}
\item $\mu$ is induced by a continuous Schauder frame over $(\Omega,\Sigma,\lambda)$ for some $\sigma$-finite positive measure $\lambda$.
\item $\mu$ admits a cyclic dilation and its operator-norm variation $\abs{\mu}$ is $\sigma$-finite.
\end{enumerate}
\end{maintheorem}

In the converse direction, the Radon--Nikod\'ym property of $X$ follows from the hypotheses; it is not an additional assumption. Cyclicity then forces an operator-valued density to have rank at most one. We also obtain associated K\"othe spaces for continuous Schauder frames while retaining their coefficient functions (Corollary~\ref{cor:associated-spaces}). This is related to the construction of shrinking and reflexive associated spaces for discrete Schauder frames by Beanland, Freeman, and Liu~\cite{BeanlandFreemanLiu2015}.

Finally, we construct compatible families of cyclic $L_p$-dilations whose inductive limits recover the underlying space and its norm (Theorem~\ref{thm:measure-inductive-dilation}). This construction applies to measure quotients whose synthesis map has a bounded linear right inverse, and we give a concrete realization for ReLU Barron spaces. Representations of Barron spaces by unions of reproducing kernel Hilbert spaces and by quotients of $L_p$-spaces are already available~\cite{EMaWu2022,BartolucciDeVitoRosascoVigogna2023,SpekHeeringaSchwenningerBrune2025}. Here the local analysis maps, synthesis maps, and multiplication projections are compatible with the same operator-valued measure, where compatibility means that the connecting maps intertwine the local analysis and synthesis operators and commute with multiplication by characteristic functions. For the Barron example, this also gives an exact variation-norm formula for the error of parameter localization while retaining the affine part.

The paper is organized as follows. In Section~2, we introduce normal parts and the localization condition and study their behavior under restrictions and quotients. Section~3 establishes James-type characterizations of shrinking, bounded completeness, and reflexivity for projection-valued measures. In Section~4, we construct cyclic dilations and give criteria for reflexive dilations in terms of weak compactness and local variation. Section~5 concerns continuous Schauder frames, frame representations of operator-valued measures, and associated K\"othe spaces. In Section~6, we study compatible cyclic dilations and their inductive limits, with applications to ReLU Barron spaces. Finally, Section~7 gives concrete examples illustrating the results and the necessity of the assumptions.

\textbf{Notation.} Banach spaces and scalars are real unless otherwise stated, and $X\neq\{0\}$. We use the following notation.
\begin{itemize}[leftmargin=2em,itemsep=2pt,parsep=0pt,topsep=3pt,partopsep=0pt]
\item $\mathcal{B}(X)$ --- bounded linear operators on $X$, with the operator norm.
\item $B_X,I_X$ --- the closed unit ball and identity operator of $X$.
\item $J_X\colon X\to X^{**}$ --- the canonical embedding.
\item $Y\hookrightarrow X$ --- $X$ contains a closed subspace isomorphic to $Y$.
\item $\mathbbm{1}_E$ --- the characteristic function of $E$.
\item $\operatorname{span} A,\overline{\operatorname{span}} A$ --- the linear span and norm-closed linear span of $A$.
\item $\overline{\operatorname{aconv}}A$ --- the norm-closed absolutely convex hull of $A$.
\item $B_1(X)$ --- weak-star limits of sequences in $J_XX$.
\item $x\otimes f$ --- $(x\otimes f)y=f(y)x$, where $x,y\in X$ and $f\in X^*$.
\end{itemize}

\section{Normal parts and localization}\label{sec:normal-localization}

Let $(\Omega,\Sigma)$ be an arbitrary measurable space.

\begin{definition}
An \emph{operator-valued measure} (OVM) on $X$ is a map $\mu\colon\Sigma\to\cB(X)$ such that $\mu(\varnothing)=0$ and
\[
\mu\Big(\bigcup_{n=1}^{\infty}E_n\Big)x=\sum_{n=1}^{\infty}\mu(E_n)x\quad(x\in X)
\]
in norm whenever $(E_n)$ is pairwise disjoint. We call this property \emph{strong countable additivity}. The measure is \emph{normalized} if $\mu(\Omega)=I_X$. A normalized OVM is a \emph{projection-valued measure} (PVM) if $\mu(E)\mu(F)=\mu(E\cap F)$ for all $E,F\in\Sigma$. For a PVM $P$, let $X(x)=\cspan\{P(E)x:E\in\Sigma\}$. The representation is \emph{cyclic} if $X=X(x)$ for some $x\in X$, and \emph{finitely generated} if $X=\overline{X(x_1)+\cdots+X(x_r)}$ for finitely many vectors $x_1,\ldots,x_r$.
\end{definition}

\begin{remark}\label{rem:ovm-bounds}
By the Orlicz--Pettis theorem, weak countable additivity implies strong countable additivity. All disjoint series defining $\mu$ converge unconditionally on each vector, and $M_\mu:=\sup_{E\in\Sigma}\norm{\mu(E)}<\infty$. For disjoint $E_1,\ldots,E_m$ and  $\abs{a_j}\leq1$,
\begin{equation}\label{eq:multiplier-bound}
\norm{\sum_{j=1}^m a_j\mu(E_j)}\leq2M_\mu.
\end{equation}
Thus $\pi(a)=\int_\Omega a\,\rmd P$ defines a bounded unital representation of the bounded measurable functions, with $\norm{\pi(a)}\leq2M_P\norm{a}_\infty$.
\end{remark}

For an OVM $\mu$, define its \emph{normal part} in the bidual by
\[
\mathcal N_\mu=\{\xi\in X^{**}:E\mapsto\xi(\mu(E)^*f)\text{ is countably additive for every }f\in X^*\}.
\]
It records which bidual vectors preserve countable additivity of the scalar coefficient measures.
\begingroup
For $\xi\in X^{**}$ and $f\in X^*$, let $\lambda^\mu_{\xi,f}(E)=\inp{\xi}{\mu(E)^*f}$. By~\eqref{eq:multiplier-bound}, $\lambda^\mu_{\xi,f}$ is a bounded finitely additive scalar measure. Let $r$ denote the norm-one projection onto the countably additive part in the Yosida--Hewitt decomposition~\cite{YosidaHewitt1952}. It commutes with restriction to measurable sets. Define $Q_\mu\colon X^{**}\to X^{**}$ by
\begin{equation}\label{eq:ovm-normal-projection}
\inp{Q_\mu\xi}{f}=(r\lambda^\mu_{\xi,f})(\Omega).
\end{equation}
Then $\norm{Q_\mu}\leq2M_\mu$.
\par\endgroup

For a closed subspace $Y\subseteq X$, let $j_Y\colon Y\to X$ be inclusion and let $q\colon X\to X/Y$ be the quotient map. For a PVM $P$, we call $Y$ $P$-invariant if $P(E)Y\subseteq Y$ for every $E\in\Sigma$. In this case, let $P_Y(E)=P(E)|_Y$ and let $\bar P$ be the induced PVM on $X/Y$.

\begin{lemma}\label{lem:normal-restriction}
For every PVM $P$, the operator $Q_P$ is a bounded projection onto $\mathcal N_P$, and $B_1(X)\subseteq\mathcal N_P$.
\begingroup
If a bounded operator $A\colon X\to Z$ and a PVM $R$ satisfy $AP(E)=R(E)A$ for every $E\in\Sigma$, then we have $ Q_RA^{**}=A^{**}Q_P $.
\endgroup
If $Y$ is closed and invariant, then
\begin{equation}\label{eq:normal-restriction}
\mathcal N_P\cap j_Y^{**}Y^{**}=j_Y^{**}\mathcal N_{P_Y}.
\end{equation}
\begingroup
Moreover, we have the short exact sequence
\[
0\longrightarrow\mathcal N_{P_Y}\xrightarrow{j_Y^{**}}\mathcal N_P\xrightarrow{q^{**}}\mathcal N_{\bar P}\longrightarrow0.
\]
\endgroup
For every closed subspace $Y\subseteq X$,
\begin{equation}\label{eq:baire-hereditary}
B_1(X)\cap j_Y^{**}Y^{**}=j_Y^{**}B_1(Y).
\end{equation}
\end{lemma}

\begin{proof}
\begingroup
Multiplicativity gives $\lambda^P_{\xi,P(E)^*f}=\lambda^P_{\xi,f}|_E$. Restriction commutes with $r$, so $ \inp{Q_P\xi}{P(E)^*f}=(r\lambda^P_{\xi,f})(E)$. Thus $Q_P$ maps into $\mathcal N_P$ and fixes every element of $\mathcal N_P$. The coefficient identity $\lambda^R_{A^{**}\xi,g}=\lambda^P_{\xi,A^*g}$ proves $ Q_RA^{**}=A^{**}Q_P $.
\par\endgroup

Suppose that $J_Xx_n\xrightarrow{w^*}\xi$. For each $f\in X^*$, the measures $\nu_n(E)=f(P(E)x_n)$ converge at every $E$ to $\xi(P(E)^*f)$. Let
\[
\lambda=\sum_{n=1}^{\infty}2^{-n}\frac{\abs{\nu_n}}{1+\abs{\nu_n}(\Omega)}.
\]
Then $\lambda(\Omega)\leq1$ and $\abs{\nu_n}\ll\lambda$ for every $n$. By the Vitali--Hahn--Saks theorem \cite[Chapter~I]{DiestelUhl1977}, the limit is countably additive. Thus $\xi\in \mathcal N_P$. Since every functional on $Y$ extends to $X$, the  identity $P(E)j_Y=j_YP_Y(E)$ proves~\eqref{eq:normal-restriction}.
\begingroup
The operator $j_Y^{**}$ is injective. By~\eqref{eq:normal-restriction} and $\ker q^{**}=j_Y^{**}Y^{**}$, we have
\[
\ker(q^{**}|_{\mathcal N_P})=\mathcal N_P\cap j_Y^{**}Y^{**}=j_Y^{**}\mathcal N_{P_Y}.
\]
For $\eta\in\mathcal N_{\bar P}$, choose $\xi\in X^{**}$ with $q^{**}\xi=\eta$. By the intertwining identity, $q^{**}Q_P\xi=Q_{\bar P}\eta=\eta$, so $q^{**}|_{\mathcal N_P}$ is onto $\mathcal N_{\bar P}$.
\par\endgroup

One inclusion in~\eqref{eq:baire-hereditary} follows from weak-star continuity of $j_Y^{**}$. For the other, let $J_Xx_n\xrightarrow{w^*}j_Y^{**}\eta$. Then $qx_n\to0$ weakly. By Mazur's theorem, for every $k$ there is a finite convex combination $v_k$ of $\{x_n:n\geq k\}$ such that $\norm{qv_k}<1/k$. Choose $y_k\in Y$ with $\norm{v_k-y_k}<2/k$. Since $J_Xv_k\xrightarrow{w^*}j_Y^{**}\eta$, we have $J_Yy_k\xrightarrow{w^*}\eta$.
\end{proof}

For a countably additive vector measure $m\colon\Sigma\to X$, a set $E\in\Sigma$ is \emph{$m$-null} if $m(F)=0$ for every measurable $F\subseteq E$. A finite positive measure $\lambda$ is a \emph{control measure} for $m$ if its null sets are precisely the $m$-null sets.

For cyclic representations, normal bidual vectors can be localized in the original space.

\begin{proposition}\label{prop:cyclic-normal}
Let $P$ be cyclic on $X$, with cyclic vector $u$. There are a finite positive measure $\lambda$ and $\varphi\in X^*$ such that
\begin{equation}\label{eq:faithful-control}
\lambda(E)=\varphi(P(E)u)\quad\text{and}\quad \lambda(E)=0\ \Longleftrightarrow\ P(E)=0.
\end{equation}
For every $\xi\in \mathcal N_P$, there are increasing $E_n\in\Sigma$ with $\bigcup_nE_n=\Omega$ and $P(E_n)^{**}\xi\in J_XX$. In particular, $\mathcal N_P=B_1(X)$.
\end{proposition}

\begin{proof}
We first choose a control measure and express the coefficient measures of each normal vector through one scalar density. We then truncate this density to obtain vectors in $X$. Define $m(E)=P(E)u$. By Rybakov's theorem~\cite{Rybakov1970}, choose $f_0\in X^*$ such that $\lambda=\abs{f_0m}$ is a control measure for $m$. If $\lambda(E)=0$, then $P(E)P(F)u=0$ for every $F\in\Sigma$, so $P(E)=0$ by cyclicity. The converse is immediate. Since every coefficient measure of $P$ is absolutely continuous with respect to $\lambda$, the representation $\pi$ from Remark~\ref{rem:ovm-bounds} satisfies $\pi(a)=\pi(b)$ whenever bounded measurable functions $a,b$ agree $\lambda$-a.e. Hence $\pi$ induces a bounded unital representation of $L_\infty(\lambda)$, still denoted by $\pi$. Let $s=\rmd(f_0m)/\rmd\lambda$, so that $\abs{s}=1$ a.e., and define $\varphi=\pi(s)^*f_0$. Then $\varphi(P(E)u)=\int_E s^2\,\rmd\lambda=\lambda(E)$.

For $f\in X^*$, let $g_f=\rmd(fm)/\rmd\lambda$. Suppose that $g_f$ is bounded on $B\in\Sigma$. For every bounded measurable $a$, we have
\begin{align*}
f(P(B)\pi(a)u)
&=\int_B ag_f\,\rmd\lambda
=\varphi(\pi(g_f\1_B)\pi(a)u).
\end{align*}
By cyclicity, it follows that
\begin{equation}\label{eq:cyclic-dual-density}
P(B)^*f=\pi(g_f\1_B)^*\varphi.
\end{equation}

Fix $\xi\in \mathcal N_P$, and let $h$ be the Radon--Nikod\'ym derivative of the measure $E\mapsto\xi(P(E)^*\varphi)$ with respect to $\lambda$. For $B_n=\{\abs{g_f}\leq n\}$, equation~\eqref{eq:cyclic-dual-density} implies
\[
\xi(P(E\cap B_n)^*f)=\int_{E\cap B_n}hg_f\,\rmd\lambda.
\]
The variation of the measure on the left bounds $\int_{B_n}\abs{hg_f}\,\rmd\lambda$, uniformly in $n$. Thus $hg_f\in L_1(\lambda)$. Since both coefficient measures are countably additive and vanish on $\lambda$-null sets, we obtain the coefficient identity
\begin{equation}\label{eq:normal-density}
\xi(P(E)^*f)=\int_E hg_f\,\rmd\lambda\quad(E\in\Sigma,\ f\in X^*).
\end{equation}
Choose a representative of $h$ finite outside a measurable null set $A$, and let $E_n=\{\abs{h}\leq n\}\cup A$. For $x_n=\pi(h\1_{E_n})u$, the coefficient identity~\eqref{eq:normal-density} implies
\[
f(x_n)=\int_{E_n}hg_f\,\rmd\lambda=\xi(P(E_n)^*f).
\]
Thus $J_Xx_n=P(E_n)^{**}\xi$ and $\norm{x_n}=\norm{P(E_n)^{**}\xi}\leq M_P\norm{\xi}$. By countable additivity, $J_Xx_n\xrightarrow{w^*}\xi$. By Lemma~\ref{lem:normal-restriction}, $\mathcal N_P=B_1(X)$.
\end{proof}

The following description of invariant subspaces will be used to compare cyclic dilations.

\begin{lemma}\label{lem:invariant-cyclic}
Let $P$ be cyclic on $X$, with cyclic vector $u$. Every closed invariant subspace $Y$ equals $P(A)X$ for some $A\in\Sigma$, and $P(A)u$ is cyclic for $P_Y$.
\end{lemma}

\begin{proof}
Use $\lambda$ and $\varphi$ from Proposition~\ref{prop:cyclic-normal}. For $y\in Y$, let $h_y$ be the derivative of $E\mapsto\varphi(P(E)y)$, and let $A_y=\{h_y\neq0\}$. By the coefficient identity~\eqref{eq:normal-density}, applied to $J_Xy$, we have $P(A_y)y=y$. For $D_n=\{1/n\leq\abs{h_y}\leq n\}$, the same identity implies $P(D_n)y=\pi(h_y\1_{D_n})u$. Thus
\[
P(D_n)u=\pi(h_y^{-1}\1_{D_n})y\in Y.
\]
Since $Y$ is closed, $P(A_y)u\in Y$ by strong countable additivity. Thus $X(y)=P(A_y)X\subseteq Y$.

Let $ \alpha=\sup\Big\{\lambda\Big(\bigcup_{y\in F}A_y\Big):F\subseteq Y\text{ finite}\Big\}$. Choose increasing finite sets $F_n\subseteq Y$ such that $B_n=\bigcup_{y\in F_n}A_y$ satisfies $\lambda(B_n)>\alpha-1/n$. Let $A=\bigcup_nB_n$, so $\lambda(A)=\alpha$. For every $y\in Y$, continuity from below gives
\[
\lambda(A\cup A_y)=\lim_n\lambda(B_n\cup A_y)\leq\alpha=\lambda(A),
\]
and hence $\lambda(A_y\setminus A)=0$. Since $P(A_y)X\subseteq Y$, we have $P(B_n)u\in Y$, and strong countable additivity gives $P(A)u\in Y$. Thus $P(A)X\subseteq Y$, whereas $P(A)y=y$ for every $y\in Y$.
\end{proof}

\begingroup
We say that a PVM $P$ satisfies the \emph{localization property} if
\begin{equation*}\tag{L}\label{eq:localization}
\forall\xi\in\mathcal N_P,\quad \exists E_n\uparrow\Omega\quad\text{such that}\quad P(E_n)^{**}\xi\in J_XX\quad(n\geq1).
\end{equation*}
This property is motivated by the truncation formulation of bounded completeness for continuous Schauder frames in~\cite[Definition~3.6]{EisnerFreeman2020}. The next proposition identifies this property with coverage by cyclic normal parts and allows finitely many generators in place of one.
\par\endgroup

\begin{proposition}\label{prop:normal-coverage}
Let $P$ be a PVM on $X$. The following conditions are equivalent:
\begin{enumerate}
\item $\mathcal N_P=\bigcup_{x\in X}j_{X(x)}^{**}\mathcal N_{P_{X(x)}}$;
\item $\mathcal N_P=\bigcup_Yj_Y^{**}\mathcal N_{P_Y}$, where $Y$ ranges over the closed finitely generated invariant subspaces of $X$;
\item $P$ satisfies the localization property~\eqref{eq:localization}.
\end{enumerate}
These conditions imply $\mathcal N_P=B_1(X)$. They hold for every finitely generated representation. For every closed $P$-invariant subspace $Y$, they hold for $P$ if and only if they hold for both $P_Y$ and $\bar P$.
\end{proposition}

\begin{proof}
\begingroup
First suppose that $P$ satisfies (3), and let $Y$ be closed and invariant. For $\eta\in\mathcal N_{P_Y}$, apply (3) to $j_Y^{**}\eta$. Then
\[
j_Y^{**}P_Y(E_n)^{**}\eta=P(E_n)^{**}j_Y^{**}\eta\in J_XX\cap j_Y^{**}Y^{**}=J_XY,
\]
so $P_Y$ satisfies (3). For $\eta\in\mathcal N_{\bar P}$, Lemma~\ref{lem:normal-restriction} gives $\xi\in\mathcal N_P$ with $q^{**}\xi=\eta$. Applying (3) to $\xi$, we have $\bar P(E_n)^{**}\eta=q^{**}P(E_n)^{**}\xi\in J_{X/Y}(X/Y)$. Thus $\bar P$ satisfies (3).

Conversely, suppose that both $P_Y$ and $\bar P$ satisfy (3), and let $\xi\in\mathcal N_P$.
\endgroup
Since $q^{**}\xi$ is normal for the quotient representation, there is a measurable partition $(D_n)$ of $\Omega$ such that $q^{**}P(D_n)^{**}\xi\in J_{X/Y}(X/Y)$. Choose $y_n\in X/Y$ such that $J_{X/Y}y_n=q^{**}P(D_n)^{**}\xi$. Since $\bar P(D_n)^{**}J_{X/Y}y_n=J_{X/Y}y_n$, we have $\bar P(D_n)y_n=y_n$. Choose $t_n\in X$ with $qt_n=y_n$, and let $x_n=P(D_n)t_n$. Then $x_n\in P(D_n)X$, $qx_n=y_n$, and
\[
P(D_n)^{**}\xi-J_Xx_n\in\ker q^{**}=j_Y^{**}Y^{**}.
\]
Both terms on the left are normal, so~\eqref{eq:normal-restriction} shows $P(D_n)^{**}\xi-J_Xx_n=j_Y^{**}\eta_n$ for some $\eta_n\in\mathcal N_{P_Y}$. Applying (3) in $Y$, choose a measurable partition $(F_{n,k})_k$ of $\Omega$ such that $P_Y(F_{n,k})^{**}\eta_n\in J_YY$. Then
\[
P(D_n\cap F_{n,k})^{**}\xi
=J_XP(D_n\cap F_{n,k})x_n+j_Y^{**}P_Y(D_n\cap F_{n,k})^{**}\eta_n\in J_XX.
\]
The sets $D_n\cap F_{n,k}$ form a countable partition. Their finite unions verify (3) on $X$.

By Proposition~\ref{prop:cyclic-normal} and extension stability, every finitely generated representation satisfies (3) by induction. At each step, use the closed invariant subspace generated by the first $r-1$ vectors and the cyclic quotient.

Clearly, (1) implies (2). Assume (2), and let $\xi\in\mathcal N_P$. Choose $Y$ as in (2) and $\eta\in\mathcal N_{P_Y}$ with $\xi=j_Y^{**}\eta$. Since $P_Y$ satisfies (3), there are $E_n\uparrow\Omega$ with $P_Y(E_n)^{**}\eta\in J_YY$. Hence
\[
P(E_n)^{**}\xi=j_Y^{**}P_Y(E_n)^{**}\eta\in j_Y^{**}J_YY=J_XY\subseteq J_XX.
\]
Thus (2) implies (3). Conversely, assume (3), let $E_0=\varnothing$, and let $J_Xz_n=P(E_n\setminus E_{n-1})^{**}\xi$. We now choose $x$ so that every $z_n$ belongs to $X(x)$. Since $\norm{z_n}\leq M_P\norm{\xi}$, the series $x=\sum_n2^{-n}z_n$ converges in norm. Since $P(E_n\setminus E_{n-1})x=2^{-n}z_n$, every $z_n$ belongs to $X(x)$. The partial sums of $\sum_nJ_Xz_n$ converge weak-star to $\xi$, so $\xi\in j_{X(x)}^{**}X(x)^{**}$. By~\eqref{eq:normal-restriction}, $\xi\in j_{X(x)}^{**}\mathcal N_{P_{X(x)}}$. This proves (1). Finally, (3) expresses every $\xi\in \mathcal N_P$ as a weak-star limit of the bounded sequence $P(E_n)^{**}\xi$ in $J_XX$. By Lemma~\ref{lem:normal-restriction}, $\mathcal N_P=B_1(X)$.
\end{proof}

\section{James-type theorems for projection-valued measures}

The coordinate projections of an unconditional Schauder basis form a PVM.

Under the localization property~\eqref{eq:localization}, we will show that shrinking corresponds to $\mathcal N_P=X^{**}$, whereas bounded completeness corresponds to $\mathcal N_P=J_XX$. The first equality requires every bidual vector to be normal, and the second requires every normal vector to be represented in $X$.

\begin{definition}\label{def:shrinking-bounded-complete}
An OVM $\mu$ is \emph{shrinking} if $E\mapsto\mu(E)^*$ is strongly countably additive on $X^*$. We call a PVM $P$ \emph{boundedly complete} when the following condition holds: if $(E_n)$ is a pairwise disjoint sequence in $\Sigma$, $z_n\in P(E_n)X$, and $\sup_{N\geq1}\norm{\sum_{n=1}^Nz_n}<\infty$, then $\lim_{N\to\infty}\norm{\sum_{n=1}^Nz_n-z}=0$ for some $z\in X$.
\end{definition}

\begingroup
If $(e_n,e_n^*)$ is an unconditional Schauder basis of $X$, then
\[
P(A)x=\sum_{n\in A}e_n^*(x)e_n\qquad(A\subseteq\N)
\]
defines a PVM. Taking $E_n=\{1,\ldots,n\}$ verifies~\eqref{eq:localization}. The shrinking and bounded completeness properties in Definition~\ref{def:shrinking-bounded-complete} coincide with those of the basis.

On $\Omega=\{\omega\}$, set $P(\varnothing)=0$ and $P(\Omega)=I_X$. This PVM is shrinking and boundedly complete for every $X$, and $\mathcal N_P=X^{**}$. If $E_n\uparrow\Omega$, then $E_n=\Omega$ for all sufficiently large $n$. Hence the localization property~\eqref{eq:localization} holds if and only if $X$ is reflexive.
\par\endgroup

\begin{proposition}\label{prop:no-l1-general}
If $\ell_1\not\hookrightarrow X$, every OVM on $X$ is shrinking.
\end{proposition}

\begin{proof}
Fix $f\in X^*$, disjoint $(E_n)$, and $U=\bigcup_nE_n$. By~\eqref{eq:multiplier-bound}, the operator $Tx=(f(\mu(E_n)x))_n$ from $X$ to $\ell_1$ satisfies $\norm{T}\leq2M_\mu\norm{f}$. By Rosenthal's theorem~\cite{Rosenthal1974}, every bounded sequence in $X$ has a weakly Cauchy subsequence $(x_k)$. Then $T(x_j-x_k)\xrightarrow{w}0$ as $j,k\to\infty$, and the Schur property of $\ell_1$ gives $\lim_{j,k\to\infty}\norm{Tx_j-Tx_k}_1=0$. Thus $T$ is compact. If $R_N$ deletes the first $N$ coordinates, then
\[
\norm{\mu(U)^*f-\sum_{n=1}^N\mu(E_n)^*f}
\leq\sup_{x\in B_X}\sum_{n>N}\abs{f(\mu(E_n)x)}
=\norm{R_NT}\longrightarrow0.
\]
\end{proof}

\begin{remark}\label{rem:finite-subsums}
Suppose that $z_n\in P(E_n)X$, the sets $E_n$ are pairwise disjoint, and $C=\sup_N\norm{\sum_{n=1}^Nz_n}<\infty$. For finite $A$ and $N\geq\max A$, we have
\[
\norm{\sum_{n\in A}z_n}
=\norm{P\Big(\bigcup_{n\in A}E_n\Big)\sum_{n=1}^Nz_n}\leq M_PC.
\]
Thus $\sum_nz_n$ is \emph{weakly unconditionally Cauchy}, meaning that $\sum_n\abs{f(z_n)}<\infty$ for every $f\in X^*$, and here this sum is at most $2M_PC\norm{f}$. By the same estimate on tails, norm convergence implies unconditional convergence.
\end{remark}

\Needspace{19\baselineskip}
\begin{theorem}\label{thm:representation-james}
Let $P$ satisfy the localization property~\eqref{eq:localization}. Then:
\begin{enumerate}
\item The following are equivalent:
    \begin{enumerate}[label=(\roman*),ref=\theenumi(\roman*),leftmargin=2em,itemsep=0pt,parsep=0pt,topsep=3pt,partopsep=0pt,beginpenalty=10000,midpenalty=10000]
\item\label{item:james-shrinking} $P$ is shrinking.
\item\label{item:james-no-l1} $\ell_1\not\hookrightarrow X$.
\item\label{item:james-normal-full} $\mathcal N_P=X^{**}$.
    \end{enumerate}
\item The following are equivalent:
    \begin{enumerate}[label=(\roman*),ref=\theenumi(\roman*),leftmargin=2em,itemsep=0pt,parsep=0pt,topsep=3pt,partopsep=0pt,beginpenalty=10000,midpenalty=10000]
\item $P$ is boundedly complete.
\item $c_0\not\hookrightarrow X$.
\item $X$ is weakly sequentially complete.
\item $\mathcal N_P=J_XX$.
    \end{enumerate}
\item $X$ is reflexive if and only if $P$ is shrinking and boundedly complete.
\end{enumerate}
In particular, these assertions apply to cyclic and finitely generated representations.
\end{theorem}

\begin{proof}
The equivalence of~\ref{item:james-shrinking} and~\ref{item:james-normal-full} follows from the definition of $\mathcal N_P$ and the Orlicz--Pettis theorem. Suppose that $P$ is shrinking. By Proposition~\ref{prop:normal-coverage}, $B_1(X)=X^{**}$. If a closed subspace $Y\subseteq X$ were isomorphic to $\ell_1$, equation~\eqref{eq:baire-hereditary} would imply $B_1(Y)=Y^{**}$. Since $\ell_1$ is weakly sequentially complete, $B_1(Y)=J_YY$, contradicting nonreflexivity of $\ell_1$. The implication from~\ref{item:james-no-l1} to~\ref{item:james-shrinking} is Proposition~\ref{prop:no-l1-general}.

Suppose that $P$ is boundedly complete and $\xi\in \mathcal N_P$. Choose $E_n\uparrow\Omega$ as in~\eqref{eq:localization}, let $E_0=\varnothing$, and let $J_Xz_n=P(E_n\setminus E_{n-1})^{**}\xi$. Then
\[
\norm{\sum_{n=1}^Nz_n}=\norm{P(E_N)^{**}\xi}\leq M_P\norm{\xi}.
\]
By bounded completeness, $\sum_nz_n$ converges in $X$. Its weak-star sum is $\xi$, so $\xi\in J_XX$. Therefore $\mathcal N_P=J_XX$, and $B_1(X)=J_XX$ means that $X$ is weakly sequentially complete. Such a space contains no copy of $c_0$, since the partial sums of the unit vector basis of $c_0$ form a weakly Cauchy sequence with no weak limit. Conversely, if $c_0\not\hookrightarrow X$, then every series in the definition of bounded completeness is weakly unconditionally Cauchy by Remark~\ref{rem:finite-subsums}. By the Bessaga--Pe\l czy\'nski theorem~\cite{BessagaPelczynski1958}, it converges unconditionally in norm. These implications prove (2).

By (1) and (2), shrinking and bounded completeness together imply $X^{**}=\mathcal N_P=J_XX$. Conversely, a reflexive space contains neither $\ell_1$ nor $c_0$, so (1) and (2) apply.
\end{proof}

These results also admit the following local form.

\begin{corollary}\label{cor:local-representation-james}
Suppose that every countable subset of $X$ is contained in a closed invariant subspace $Y$ such that $P_Y$ satisfies the localization property~\eqref{eq:localization}. Then:
\begin{enumerate}
\item $P$ is shrinking if and only if $\ell_1\not\hookrightarrow X$.
\item $P$ is boundedly complete if and only if $c_0\not\hookrightarrow X$, equivalently if and only if $X$ is weakly sequentially complete.
\item $X$ is reflexive if and only if $P$ is shrinking and boundedly complete.
\end{enumerate}
\end{corollary}

\begin{proof}
Shrinking and bounded completeness pass to closed invariant subspaces, and for shrinking we extend functionals on $Y$ by the Hahn--Banach theorem. If $P$ is shrinking or boundedly complete, place the basis of a putative copy of $\ell_1$ or $c_0$ in one such $Y$ and apply Theorem~\ref{thm:representation-james}. The converse implications follow from Proposition~\ref{prop:no-l1-general} and Remark~\ref{rem:finite-subsums}. A weakly Cauchy sequence lies in one such $Y$, and weak sequential completeness follows from the same theorem when $c_0\not\hookrightarrow X$. If $P$ is shrinking and boundedly complete, every bounded sequence lies in a reflexive $Y$ and has a weakly convergent subsequence. By the Eberlein--\v Smulian theorem, $X$ is reflexive.
\end{proof}

A Banach lattice $E$ is \emph{order continuous} if $0\leq h_\alpha\downarrow0$ in $E$ implies $\norm{h_\alpha}_E\to0$. The following remark relates the preceding results to Banach lattices and Banach $C(K)$-modules.

\begin{remark}\label{rem:lattice-james}
\emph{Lattice realization.} For cyclic $P$, equation~\eqref{eq:faithful-control} identifies $P(\Sigma)$ with the complete measure algebra of a finite measure. Let $(P(E_\alpha))$ be an increasing net, with measurable sets considered modulo $\lambda$-null sets. Choose increasing indices $\alpha_n$ such that $\lambda(E_{\alpha_n})\to\beta:=\sup_\alpha\lambda(E_\alpha)$, and let $B=\bigcup_nE_{\alpha_n}$. By directedness, $\lambda(E_{\alpha_n}\cup E_\alpha)\leq\beta$ for every $n$ and $\alpha$, so
\[
\lambda(B\cup E_\alpha)=\lim_n\lambda(E_{\alpha_n}\cup E_\alpha)\leq\beta=\lambda(B).
\]
Thus $P(B)=\sup_\alpha P(E_\alpha)$. For $\alpha\geq\alpha_n$ and $x\in X$, we have
\[
\norm{(P(B)-P(E_\alpha))x}\leq M_P\norm{(P(B)-P(E_{\alpha_n}))x}.
\]
The right-hand side tends to zero by strong countable additivity, so the net converges strongly to $P(B)$. Thus $P(\Sigma)$ is Bade complete. Every nonzero $P(E)X$ is cyclic with vector $P(E)u$, so the Boolean algebra has uniform multiplicity one. Its classical representation therefore equips $X$ with an equivalent order-continuous Banach lattice norm, with $u$ positive and quasi-interior and $P(\Sigma)$ equal to the band projections, by~\cite[Proposition~2.4 and Lemma~2.9]{KitoverOrhon2018}.

Conversely, the geometric conclusions above hold for the band-projection representation of every order-continuous Banach lattice. Here measurable sets are replaced by elements of the complete Boolean algebra of band projections, and countable unions by joins. For any sequence $(x_n)$, let $u=\sum_n2^{-n}\abs{x_n}/(1+\norm{x_n})$, and let
\[
I_u=\{x\in X:\abs{x}\leq cu\text{ for some }c\geq0\}.
\]
Since $\abs{x_n}\leq2^n(1+\norm{x_n})u$, every $x_n$ belongs to $I_u$. An order-continuous Banach lattice has the principal projection property. By Freudenthal's spectral theorem, for every $x\in I_u$ and $\varepsilon>0$, there is a finite linear combination $s$ of vectors $Ru$, with $R$ a band projection, such that $\abs{x-s}\leq\varepsilon u$. Hence
\[
\overline{I_u}=\cspan\{Ru:R\text{ is a band projection on }X\}.
\]
This closed ideal is cyclic and contains every $x_n$, so Corollary~\ref{cor:local-representation-james} applies.

A Banach lattice is a \emph{KB-space} if every increasing norm-bounded positive sequence converges in norm. This is equivalent to weak sequential completeness by~\cite[Theorem~2.18 and the proof of Lemma~2.19]{KitoverOrhon2017}. Thus bounded completeness also identifies the KB property.

\emph{Comparison with Banach $C(K)$-modules.} Suppose that a finite positive measure $\lambda$ satisfies $\lambda(E)=0$ if and only if $P(E)=0$. Such a measure exists for finitely generated representations by adding the control measures of the generators. Then $\pi$ induces a faithful representation of $L_\infty(\lambda)\cong C(K)$, where $K$ is hyperstonian, and the preceding argument gives Bade completeness. Let $p\in C(K)^{**}$ be the central projection associated with the predual $L_1(\lambda)$, acting on $X^{**}$ through the Arens extension of $\pi$. By~\cite[Lemma~2.13(3)]{KitoverOrhon2017}, we have $\mathcal N_P=pX^{**}$. The localization property~\eqref{eq:localization} is weaker than $pX^{**}=J_XX$, and Theorem~\ref{thm:representation-james} characterizes this equality by bounded completeness under this property. For example, the coordinate PVM on $c_0$ satisfies~\eqref{eq:localization} and has $\mathcal N_P=\ell_\infty\neq J_{c_0}c_0$.

The reflexivity and weak sequential completeness characterizations for finitely generated modules and modules of finite multiplicity were established in~\cite{KitoverOrhon2014,KitoverOrhon2017}. Under the hyperstonian and Bade completeness assumptions, the equality $pX^{**}=J_XX$ already has a three-space property for closed submodules by~\cite[Lemma~2.16]{KitoverOrhon2017}. Proposition~\ref{prop:normal-coverage} concerns the weaker localization property~\eqref{eq:localization}. It identifies this property with coverage by cyclic normal parts and proves its stability under invariant subspaces, quotients, and extensions. Example~\ref{ex:vector-sequences} gives representations satisfying this property outside finite multiplicity.
\end{remark}

\begin{remark}\label{rem:normal-quotients}
\begingroup
By Lemma~\ref{lem:normal-restriction}, $Q_P$ is a bounded projection onto $\mathcal N_P$ and fixes $J_XX$ for every PVM $P$. Therefore the short exact sequence of Banach spaces
\[
0\longrightarrow\mathcal N_P/J_XX\longrightarrow X^{**}/J_XX\longrightarrow X^{**}/\mathcal N_P\longrightarrow0
\]
splits, with Banach-space isomorphisms
\[
X^{**}/J_XX\cong(\mathcal N_P/J_XX)\oplus\ker Q_P\qquad\text{and}\qquad \ker Q_P\cong X^{**}/\mathcal N_P.
\]
The first isomorphism sends $[\xi]$ to $([Q_P\xi],(I-Q_P)\xi)$, and the second is induced by the quotient map.
\par\endgroup
The first quotient consists of normal bidual vectors modulo those represented in $X$, and the last quotient records failure of normality. By the Orlicz--Pettis theorem, $P$ is shrinking if and only if $X^{**}/\mathcal N_P=\{0\}$. If $P$ satisfies the localization property~\eqref{eq:localization}, then by Theorem~\ref{thm:representation-james}, $P$ is boundedly complete if and only if $\mathcal N_P/J_XX=\{0\}$. Under this assumption, the two end terms vanish precisely when $P$ is shrinking and boundedly complete, equivalently when $X$ is reflexive.

Strong countable additivity of $P^{**}$ implies that $P$ is shrinking and boundedly complete for every PVM $P$. The first assertion follows from the Orlicz--Pettis theorem. For the second, if $z_n\in P(E_n)X$ have bounded partial sums, choose a weak-star cluster point $\xi$ of those sums and let $U=\bigcup_nE_n$. Then $P(E_n)^{**}\xi=J_Xz_n$ and $P(U)^{**}\xi=\xi$, so strong countable additivity implies norm convergence of $\sum_nJ_Xz_n$. Under the hypotheses of Theorem~\ref{thm:representation-james} or Corollary~\ref{cor:local-representation-james}, this is equivalent to reflexivity of $X$.
\end{remark}

\section[Geometric dilations of operator-valued measures]{Geometric dilations of operator-valued measures}\label{sec:cyclic-dilations}

We next transfer these results to OVMs through their dilations.

\subsection[Minimal dilations and normal parts]{Minimal dilations and normal parts}

\begin{definition}
Let $\mu\colon\Sigma\to\cB(X)$ be a normalized OVM. A \emph{dilation} of $\mu$ consists of a Banach space $Z$, a PVM $P\colon\Sigma\to\cB(Z)$, and bounded operators $T\colon X\to Z$ and $S\colon Z\to X$ such that
\[
ST=I_X \quad \text{and}\quad \mu(E)=SP(E)T\quad(E\in\Sigma).
\]
The dilation is \emph{cyclic} if $P$ has a cyclic vector.
\end{definition}

Thus $T$ is an isomorphic embedding and $TS$ is a projection onto $T(X)$. We first recall the minimal dilation of~\cite[Definition~2.28 and Theorem~2.30]{HanLarsonLiuLiu2014}. Let $\mathcal M_\mu$ be the linear span of the vector measures $\mu_{E,x}(F)=\mu(E\cap F)x$, and let $Z_\alpha$ be its completion in the norm $\norm{\eta}_\alpha=\sup_{F\in\Sigma}\norm{\eta(F)}$. The maps defined by
\[
(P_\alpha(E)\eta)(F)=\eta(E\cap F),\quad T_\alpha x=\mu_{\Omega,x},\quad\text{and}\quad S_\alpha\eta=\eta(\Omega)
\]
extend to bounded operators, with $P_\alpha$ a normalized, strongly countably additive PVM. Moreover, $S_\alpha T_\alpha=I_X$ and $S_\alpha P_\alpha(E)T_\alpha=\mu(E)$. We call this the minimal $\alpha$-dilation of $\mu$.

\begin{proposition}\label{prop:minimal-cyclic}
A normalized OVM $\mu$ admits a cyclic dilation if and only if its minimal $\alpha$-dilation $P_\alpha$ is cyclic.
\end{proposition}

\begin{proof}
One implication is immediate. Conversely, let $(Z,P,S,T)$ be a cyclic dilation and let $W=\cspan\{P(E)Tx:E\in\Sigma,\ x\in X\}$. By Lemma~\ref{lem:invariant-cyclic}, the restriction of $P$ to $W$ has a cyclic vector $v$. On the algebraic generating span, define $(Qz)(F)=SP(F)z$. Since $Q(P(E)Tx)=\mu_{E,x}$, this defines a map onto $\mathcal M_\mu$, and
\[
\norm{Qz}_\alpha=\sup_{F\in\Sigma}\norm{SP(F)z}\leq\norm{S}M_P\norm{z}.
\]
Thus $Q$ extends to $W\to Z_\alpha$, with dense range and $QP(E)|_W=P_\alpha(E)Q$. Consequently,
\[
Z_\alpha=\overline{Q(W)}=\cspan\{P_\alpha(E)Qv:E\in\Sigma\},
\]
so $Qv$ is cyclic.
\end{proof}

\begingroup
\begin{definition}\label{def:ovm-bounded-complete}
A normalized OVM $\mu$ is \emph{boundedly complete} if the PVM $P_\alpha$ in its minimal $\alpha$-dilation is boundedly complete.
\end{definition}

For a PVM $P$ on $X$, define $(Ux)(E)=P(E)x$. Then we have $\mathcal M_P=U(X)$, $\norm{x}\leq\norm{Ux}_\alpha\leq M_P\norm{x}$, and $P_\alpha(E)U=UP(E). $ Thus $U\colon X\to Z_\alpha$ is an onto isomorphism, and the two notions of bounded completeness agree by Definition~\ref{def:shrinking-bounded-complete}. Shrinking need not pass from $\mu$ to $P_\alpha$, as Remark~\ref{rem:scalar-minimal} shows.

\begin{proposition}\label{prop:minimal-bounded-complete}
Let $\mu$ be a normalized OVM on $X$. If $c_0\not\hookrightarrow X$, then $\mu$ is boundedly complete.
\end{proposition}

\begin{proof}
Let $z_n\in P_\alpha(E_n)Z_\alpha$, where $(E_n)$ is pairwise disjoint, and suppose that $C=\sup_N\norm{\sum_{n=1}^Nz_n}_\alpha<\infty$. If the series does not converge, there are successive finite intervals $I_k\subseteq\N$ such that $w_k=\sum_{n\in I_k}z_n$ satisfies $\norm{w_k}_\alpha>\varepsilon$ for some $\varepsilon>0$. Let $F_k=\bigcup_{n\in I_k}E_n$, and choose $A_k\subseteq F_k$ with $\norm{w_k(A_k)}>\varepsilon$. For every finite $H\subseteq\N$, with $N\geq\max\bigcup_{k\in H}I_k$, we have
\[
\norm{\sum_{k\in H}w_k(A_k)}
=\norm{\left(\sum_{n=1}^Nz_n\right)\left(\bigcup_{k\in H}A_k\right)}
\leq C.
\]
Thus $\sum_kw_k(A_k)$ is weakly unconditionally Cauchy in $X$. The Bessaga--Pe\l czy\'nski theorem~\cite{BessagaPelczynski1958} gives unconditional convergence, contrary to $\norm{w_k(A_k)}>\varepsilon$.
\end{proof}
\endgroup

We now apply the operators $Q_\mu$ from~\eqref{eq:ovm-normal-projection} to dilations.

\begingroup
\begin{proposition}\label{prop:ovm-normal}
Let $\mu$ be a normalized OVM on $X$ with a dilation $(Z,P,S,T)$, where $P$ satisfies the localization property~\eqref{eq:localization}. Then:
\begin{enumerate}
\item $\mu$ is shrinking if and only if $\ell_1\not\hookrightarrow X$.
\item $\mu$ is boundedly complete if and only if $c_0\not\hookrightarrow X$.
\item $X$ is reflexive if and only if $\mu$ is shrinking and boundedly complete.
\end{enumerate}
In particular, the conclusions apply to finitely generated dilations.
\end{proposition}

\begin{proof}
By Lemma~\ref{lem:normal-restriction}, $Q_P$ is a projection onto $\mathcal N_P$, where $\mathcal N_P=B_1(Z)$ by Proposition~\ref{prop:normal-coverage}. The coefficient identity $\lambda^\mu_{\xi,f}=\lambda^P_{T^{**}\xi,S^*f}$ yields $Q_\mu=S^{**}Q_PT^{**}$. Since second adjoints preserve weak-star sequential limits and $ST=I_X$, we have
\[
S^{**}B_1(Z)\subseteq B_1(X)=S^{**}T^{**}B_1(X)\subseteq S^{**}\mathcal N_P.
\]
Hence $S^{**}\mathcal N_P=B_1(X)$ and $Q_\mu(X^{**})\subseteq B_1(X)$. If $\xi\in B_1(X)$, then $T^{**}\xi\in B_1(Z)\subseteq\mathcal N_P$ by Lemma~\ref{lem:normal-restriction}, so
\[
Q_\mu\xi=S^{**}Q_PT^{**}\xi=S^{**}T^{**}\xi=\xi.
\]
Thus $Q_\mu$ is a projection onto $B_1(X)$. The coefficient identity also gives $B_1(X)\subseteq\mathcal N_\mu$. Conversely, $Q_\mu\xi=\xi$ for $\xi\in\mathcal N_\mu$, because $r$ fixes its coefficient measures and $\mu(\Omega)=I_X$. Therefore $\mathcal N_\mu=B_1(X)$.

By the Orlicz--Pettis theorem, $\mu$ is shrinking if and only if $\mathcal N_\mu=X^{**}$. If $\mu$ is shrinking, the argument in the proof of Theorem~\ref{thm:representation-james} gives $\ell_1\not\hookrightarrow X$. The converse follows from Proposition~\ref{prop:no-l1-general}.

If $c_0\not\hookrightarrow X$, Proposition~\ref{prop:minimal-bounded-complete} shows that $\mu$ is boundedly complete. Conversely, suppose that $\mu$ is boundedly complete, and take $\xi\in B_1(X)$ and a sequence $(x_k)$ with $J_Xx_k\xrightarrow{w^*}\xi$. By Lemma~\ref{lem:normal-restriction}, $T^{**}\xi\in B_1(Z)\subseteq\mathcal N_P$. Proposition~\ref{prop:normal-coverage} gives $A_n\uparrow\Omega$ and $v_n\in Z$ such that
\[
J_Zv_n=P(A_n)^{**}T^{**}\xi \qquad\text{and}\qquad \norm{v_n}\leq M_P\norm{T^{**}\xi}.
\]
Let $W=\cspan\{P(E)Tx:E\in\Sigma,\ x\in X\}$. As in the proof of Proposition~\ref{prop:minimal-cyclic}, the formula $(Qz)(E)=SP(E)z$ defines a bounded map $Q\colon W\to Z_\alpha$ satisfying
\[
QT=T_\alpha\quad\text{and}\quad QP(E)|_W=P_\alpha(E)Q.
\]
For each $n$, we have $P(A_n)Tx_k\xrightarrow{w}v_n$ as $k\to\infty$, so $v_n\in W$ and $ J_{Z_\alpha}Qv_n=P_\alpha(A_n)^{**}T_\alpha^{**}\xi $. Let $A_0=\varnothing$, $v_0=0$, and $z_n=Q(v_n-v_{n-1})$. Then $z_n\in P_\alpha(A_n\setminus A_{n-1})Z_\alpha$ and
\[
\sup_N\norm{\sum_{n=1}^Nz_n}_\alpha=\sup_N\norm{Qv_N}_\alpha\leq\norm{Q}M_P\norm{T^{**}\xi}.
\]
Bounded completeness gives $Qv_n\to z$ in $Z_\alpha$ for some $z\in Z_\alpha$. Since $T_\alpha^{**}\xi\in B_1(Z_\alpha)\subseteq\mathcal N_{P_\alpha}$, normality gives $J_{Z_\alpha}Qv_n\xrightarrow{w^*}T_\alpha^{**}\xi$. Hence $T_\alpha^{**}\xi=J_{Z_\alpha}z$ and $\xi=S_\alpha^{**}T_\alpha^{**}\xi=J_XS_\alpha z$. Thus $B_1(X)=J_XX$, so $X$ is weakly sequentially complete and contains no copy of $c_0$. This proves (2). If $\mu$ is shrinking and boundedly complete, then $X^{**}=\mathcal N_\mu=B_1(X)=J_XX$. Conversely, reflexivity excludes both $\ell_1$ and $c_0$, so (1) and (2) give (3).
\end{proof}
\endgroup

\begingroup
\begin{remark}[Normal parts under dilation]\label{rem:ovm-normal-transfer}
Under the hypotheses of Proposition~\ref{prop:ovm-normal}, its proof gives
\[
\mathcal N_\mu=B_1(X)=S^{**}\mathcal N_P,\qquad Q_\mu=S^{**}Q_PT^{**}.
\]
The operator $Q_\mu$ is a bounded projection onto $\mathcal N_\mu$.
\end{remark}
\par\endgroup

\subsection[Vector-valued multiplication representations]{Vector-valued multiplication representations}

We now compute the normal part of a vector-valued multiplication representation and determine when it satisfies~\eqref{eq:localization}.

\begingroup
A \emph{K\"othe function space} on $(\Omega,\Sigma,\lambda)$ is a Banach space $E$ of measurable scalar functions modulo equality a.e., containing a strictly positive function and satisfying
\[
\abs{h}\leq\abs{k},\quad k\in E\quad\Longrightarrow\quad h\in E,\quad\norm{h}_E\leq\norm{k}_E.
\]
Its \emph{K\"othe dual} $E^\times$ consists of the measurable scalar functions $g$, modulo equality a.e., for which
\[
\norm{g}_{E^\times}=\sup_{\norm{f}_E\leq1}\int_\Omega\abs{fg}\,\rmd\lambda<\infty,
\]
and the \emph{K\"othe bidual} is $E^{\times\times}=(E^\times)^\times$. For a Banach space $Y$, the space $E(Y)$ consists of the strongly measurable functions $h\colon\Omega\to Y$ with $\norm{h(\cdot)}_Y\in E$, again modulo equality a.e., and has norm $\norm{h}_{E(Y)}=\norm{\norm{h(\cdot)}_Y}_E$. We pair $E(Y)$ and $E^{\times\times}(Y)$ with $E^\times(Y^*)$ by
\[
\inp{h}{k}=\int_\Omega\inp{h(t)}{k(t)}\,\rmd\lambda(t).
\]

A Banach space $Y$ has the \emph{Radon--Nikod\'ym property} if every countably additive $Y$-valued measure of bounded variation that is absolutely continuous with respect to a finite positive measure has a Bochner integrable density, as recalled in~\cite[Chapter~III]{DiestelUhl1977}.

\begin{proposition}\label{prop:bochner-normal}
Let $E$ be an order-continuous K\"othe function space with continuous inclusions $L_\infty(\lambda)\subseteq E\subseteq L_1(\lambda)$. Let $Y$ be a reflexive Banach space. On $X=E(Y)$, define $P(A)h=\1_Ah$. Under the integral identification of $E^{\times\times}(Y)$ with a subspace of $X^{**}$, we have
\[
\mathcal N_P=E^{\times\times}(Y).
\]
In particular, every normal bidual vector belongs to the bidual of a cyclic subrepresentation.
\end{proposition}

\begin{proof}
Passing to the completion does not change the function spaces or the multiplication projections, so we may assume that $\lambda$ is complete. By order continuity, $P$ is strongly countably additive.
\begingroup
Since $Y^*$ is reflexive, it has the Radon--Nikod\'ym property. By K\"othe--Bochner duality~\cite[proof of Corollary~21]{Hardtke2013}, we have $X^*=E^\times(Y^*)$ isometrically under the integral pairing.

Let $v\colon\Omega\to Y$ be strongly measurable. After modifying $v$ on a null set, we may assume that its range lies in a separable subspace $Y_v\subseteq Y$. Choose $(d_j^*)\subseteq B_{Y^*}$ such that $\sup_jd_j^*(y)=\norm{y}$ for every $y\in Y_v$. For $0<\varepsilon<1$, define
\par\endgroup
\[
j_{v,\varepsilon}(t)=\min\{j:d_j^*(v(t))\geq(1-\varepsilon)\norm{v(t)}\}\quad\text{and}\quad s_{v,\varepsilon}(t)=d_{j_{v,\varepsilon}(t)}^*.
\]
The level sets of $j_{v,\varepsilon}$ are measurable, so $s_{v,\varepsilon}$ is countably valued and strongly measurable.

For $v\in E^{\times\times}(Y)$, H\"older's inequality shows that $\xi_v(k)=\int_\Omega\inp{v(t)}{k(t)}\,\rmd\lambda(t)$ defines an element of $X^{**}$ with $\norm{\xi_v}\leq\norm{v}_{E^{\times\times}(Y)}$. For $0\leq b\in B_{E^\times}$, we have
\[
(1-\varepsilon)\int_\Omega b\norm{v}\,\rmd\lambda\leq\xi_v(bs_{v,\varepsilon})\leq\norm{\xi_v}.
\]
Taking the supremum over $b$ and then letting $\varepsilon\downarrow0$ proves equality of the norms. Moreover, $A\mapsto\xi_v(P(A)^*k)=\int_A\inp{v(t)}{k(t)}\,\rmd\lambda(t)$ is countably additive for every $k\in X^*$, so $\xi_v$ is normal.

Conversely, let $\xi\in \mathcal N_P$. For $A\in\Sigma$, define $m(A)\in Y^{**}=Y$ by $\inp{m(A)}{y^*}=\xi(\1_Ay^*)$. By normality, $m$ is weakly countably additive. By the Orlicz--Pettis theorem, it is norm countably additive. For a finite disjoint family $(A_j)$, choose $y_j^*\in B_{Y^*}$ which norm the vectors $m(A_j)$. Then
\begin{align*}
\sum_j\norm{m(A_j)}
&=\xi\left(\sum_j\1_{A_j}y_j^*\right)\\
&\leq\norm{\xi}\left\|\sum_j\1_{A_j}y_j^*\right\|_{E^\times(Y^*)}
 \leq\norm{\xi}\norm{\1_\Omega}_{E^\times}.
\end{align*}
Thus $m$ has bounded variation and $m\ll\lambda$. Since $Y$ has the Radon--Nikod\'ym property, there is $h\in L_1(\lambda;Y)$ such that $m(A)=\int_Ah\,\rmd\lambda$.

For every bounded strongly measurable $k:\Omega\to Y^*$, we have
\begin{equation}\label{eq:bochner-normal-pairing}
\xi(k)=\int_\Omega\inp{h(t)}{k(t)}\,\rmd\lambda(t).
\end{equation}
Indeed, normality proves this identity for bounded countably valued functions by summing over their level sets. Such functions approximate $k$ uniformly, so the identity follows by continuity. Let $0<\varepsilon<1$, and let $s=s_{h,\varepsilon}$ be the strongly measurable function constructed above. Then $\inp{h(t)}{s(t)}\geq(1-\varepsilon)\norm{h(t)}$. For $b\in E^\times_+$, let $A_n=\{b\leq n\}$. By~\eqref{eq:bochner-normal-pairing},
\begin{align*}
(1-\varepsilon)\int_{A_n}b\norm{h}\,\rmd\lambda
&\leq \xi(\1_{A_n}bs)\\
&\leq\norm{\xi}\norm{\1_{A_n}bs}_{E^\times(Y^*)}
 \leq\norm{\xi}\norm{b}_{E^\times}.
\end{align*}
By monotone convergence and then letting $\varepsilon\downarrow0$, we obtain $\norm{\norm{h}_Y}_{E^{\times\times}}\leq\norm{\xi}$. For arbitrary $k\in E^\times(Y^*)$, apply~\eqref{eq:bochner-normal-pairing} on $\{\norm{k}\leq n\}$. By normality on the left and dominated convergence on the right, the same identity holds for $k$. Hence $\xi$ is the functional represented by $h$.

Finally, on $B_n=\{\norm{h}\leq n\}$, the function $\1_{B_n}h$ belongs to $E(Y)$. These sets increase to $\Omega$ modulo a null set. By Proposition~\ref{prop:normal-coverage}, $\xi$ belongs to the bidual of a cyclic subrepresentation.
\end{proof}

\begin{remark}\label{rem:bochner-localization-dilation}
Let $E$ be as in Proposition~\ref{prop:bochner-normal} and $Y$ be a Banach space. On $Z=E(Y)$, define $P(A)h=\1_Ah$. Then $P$ satisfies the localization property~\eqref{eq:localization} if and only if $Y$ is reflexive. Sufficiency follows from Proposition~\ref{prop:bochner-normal}. For necessity, define bounded operators $C\colon Y\to Z$ and $V\colon Z\to Y$ by
\[
Cy=\1_\Omega y\quad\text{and}\quad Vh=\int_\Omega h\,\rmd\lambda.
\]
Then $VP(A)C=\lambda(A)I_Y$. For $\eta\in Y^{**}$, $f\in Z^*$, and $D_n\downarrow\varnothing$, order continuity gives
\[
\abs{\inp{P(D_n)^{**}C^{**}\eta}{f}}\leq\norm{\eta}\norm{f}\norm{\1_{D_n}}_E\longrightarrow0.
\]
Thus $C^{**}\eta\in\mathcal N_P$. By~\eqref{eq:localization}, there are $A_n\uparrow\Omega$ such that $P(A_n)^{**}C^{**}\eta\in J_ZZ$, and hence $ \lambda(A_n)\eta=V^{**}P(A_n)^{**}C^{**}\eta\in J_YY $. Since $\lambda(A_n)>0$ for some $n$, we have $\eta\in J_YY$, so $Y$ is reflexive.

If $Y$ is reflexive, let $T\colon X\to Z$ and $S\colon Z\to X$ be bounded operators with $ST=I_X$, and define $\mu(A)=SP(A)T$. By Proposition~\ref{prop:bochner-normal} and Remark~\ref{rem:ovm-normal-transfer},
\[
\mathcal N_\mu=B_1(X)=S^{**}\bigl(E^{\times\times}(Y)\bigr),
\]
where $E^{\times\times}(Y)$ is identified with a subspace of $Z^{**}$ by the integral pairing.
\end{remark}
\endgroup

\subsection[Cyclic dilations and their geometry]{Cyclic dilations and their geometry}

We now realize cyclic dilations on spaces of integrable functions with respect to vector measures.

Let $m\colon\Sigma\to X$ be a countably additive vector measure. A measurable scalar function $h$ is $m$-integrable if it is integrable with respect to $\abs{gm}$ for every $g\in X^*$, and for every $E\in\Sigma$ there is a vector $\int_Eh\,\rmd m\in X$ satisfying $g(\int_Eh\,\rmd m)=\int_Eh\,\rmd(gm)$ for all $g\in X^*$. Modulo $m$-null sets, these functions form the K\"othe function space $L_1(m)$, with norm
\[
\norm{h}_{L_1(m)}=\sup_{g\in B_{X^*}}\int_\Omega\abs{h}\,\rmd\abs{gm}.
\]
Its norm is order continuous, simple functions are dense, and the integration operator $I_mh=\int_\Omega h\,\rmd m$ is bounded. The multiplication measure $Q(E)h=\1_Eh$ is cyclic on $L_1(m)$, with vector $\1_\Omega$.

\begingroup
\begin{proposition}\label{prop:vector-measure-dilation}
Let $(Z,P,S,T)$ be a cyclic dilation of $\mu$, with cyclic vector $u$, and define $m(E)=SP(E)u$. There is a bounded operator $q\colon Z\to L_1(m)$ with dense range and $qu=\1_\Omega$. With $T_0=qT$, the following diagram commutes for every $E\in\Sigma$.
\[
\begin{CD}
    X @>{T}>> Z @>{P(E)}>> Z @>{S}>> X\\
    @| @V{q}VV @VV{q}V @|\\
    X @>{T_0}>> L_1(m) @>{Q(E)}>> L_1(m) @>{I_m}>> X
\end{CD}
\]
Consequently, $(L_1(m),Q,I_m,T_0)$ is a cyclic dilation of $\mu$. In the minimal $\alpha$-dilation, $q$ in the diagram is a lattice isomorphism for the lattice structure determined by $u$, and $U=q^{-1}$ satisfies
\[
\frac{1}{2}\norm{f}_{L_1(m)}\leq\norm{Uf}_\alpha\leq\norm{f}_{L_1(m)}\quad(f\in L_1(m)).
\]
\end{proposition}

\begin{proof}
The vector measure $m$ is countably additive. For a simple function $f$, let $w_f=(\int_\Omega f\,\rmd P)u$. Scalar total variation and~\eqref{eq:multiplier-bound} show that
\[
\norm{f}_{L_1(m)}
=\sup_{\substack{a\text{ simple}\\\abs{a}\leq1}}\norm{S\Big(\int_\Omega a\,\rmd P\Big)w_f}
\leq2M_P\norm{S}\norm{w_f}.
\]
Thus $q(w_f)=[f]_m$ is well defined and extends to a bounded operator on $Z$. Its range contains the simple functions and is dense in $L_1(m)$. The identities $qu=\1_\Omega$, $qP(E)=Q(E)q$ and $I_mq=S$ hold on the algebraic cyclic span and extend by continuity. Since $T_0=qT$, the diagram commutes and the assertions for $T_0$ follow from the dilation identities.

In the minimal dilation, the defining norm gives $ \norm{z}_\alpha=\sup_{E\in\Sigma}\norm{S_\alpha P_\alpha(E)z}$ for $z\in Z_\alpha$. This holds on $\mathcal M_\mu$ and extends by continuity, since the operators $S_\alpha$ and $P_\alpha(E)$ are contractions. Hence, for simple $f$,
\[
\norm{w_f}_\alpha=\sup_{E\in\Sigma}\left\|\int_Ef\,\rmd m\right\|
\leq\norm{f}_{L_1(m)}\leq2\norm{w_f}_\alpha.
\]
By density, $q$ is bounded below, so its range is closed and equals $L_1(m)$. The estimates for $U=q^{-1}$ follow. By Remark~\ref{rem:lattice-james}, the positive cone of $Z_\alpha$ is the closure of the vectors $U f$ with $f$ nonnegative and simple. These functions are dense in $L_1(m)_+$, so continuity of $U$ and $U^{-1}$ shows that $U$ is a lattice isomorphism.
\end{proof}

\endgroup

\begingroup
\begin{corollary}\label{cor:review-minimal-kb}
Let $\mu$ admit a cyclic dilation. Its minimal dilation space $Z_\alpha$, with the lattice structure and equivalent norm of Remark~\ref{rem:lattice-james}, is a KB-space if and only if $c_0\not\hookrightarrow X$ if and only if $\mu$ is boundedly complete.
\end{corollary}

\begin{proof}
Suppose that $c_0\not\hookrightarrow X$. By~\cite[Corollary~4.4(iv)]{OkadaRickerSanchezPerez2014}, we have $c_0\not\hookrightarrow L_1(m)$ for every countably additive vector measure $m\colon\Sigma\to X$. Since $L_1(m)$ is an order-continuous Banach lattice, it is a KB-space by Theorem~\ref{thm:representation-james} and Remark~\ref{rem:lattice-james}. By Proposition~\ref{prop:vector-measure-dilation}, $Z_\alpha$ is lattice isomorphic to such an $L_1(m)$, so it is a KB-space. Conversely, the KB property excludes $c_0$ from $Z_\alpha$ and hence from its subspace $T_\alpha X$. By Remark~\ref{rem:lattice-james} and Definition~\ref{def:ovm-bounded-complete}, $Z_\alpha$ is a KB-space if and only if $\mu$ is boundedly complete.
\end{proof}
\endgroup

The next lemma will be used to construct shrinking and reflexive cyclic dilations.

A bounded subset of a Banach space is \emph{weakly precompact} if every sequence has a weakly Cauchy subsequence. A subset $A$ of a Banach lattice is \emph{solid} if $z\in A$ and $\abs{h}\leq\abs{z}$ imply $h\in A$ for every $h\in Z$. Its \emph{solid hull} is $\{h:\abs{h}\leq\abs{z}\text{ for some }z\in A\}$.

\begin{lemma}\label{lem:multiplier-hull}
Let $P$ be a cyclic PVM on $Z$, and let $K\subseteq Z$ be weakly precompact. Then the set
\[
W=\overline{\operatorname{aconv}}\{\pi(a)z:\norm{a}_\infty\leq1,\ z\in K\}
\]
is weakly precompact. If $P$ is boundedly complete, then $W$ is weakly compact.
\end{lemma}

\begin{proof}
Let $u$ be a cyclic vector and define $m(E)=P(E)u$. Apply Proposition~\ref{prop:vector-measure-dilation} to the identity dilation $(Z,P,I_Z,I_Z)$. The resulting map $q\colon Z\to L_1(m)$ has dense range and satisfies $I_mq=I_Z$, so
\[
\norm{z}\leq\norm{I_m}\norm{qz}_{L_1(m)}\quad(z\in Z).
\]
Thus $q$ has closed range and is an isomorphism onto $L_1(m)$. Under this identification, $qP(E)=Q(E)q$ makes $\pi(a)$ multiplication by $a$. An order-continuous Banach lattice contains no copy of $C[0,1]$, so the closed convex solid hull of a weakly precompact set remains weakly precompact, as recalled in~\cite[proof of Theorem~4.5]{AvilesGuiraoLajaraRodriguezTradacete2016}. This applies to $W$, since the multiplier orbit of $K$ is its solid hull. If $P$ is boundedly complete, then $Z$ is weakly sequentially complete by Theorem~\ref{thm:representation-james}. The Eberlein--\v Smulian theorem therefore shows that $W$ is weakly compact.
\end{proof}

\begingroup
\begin{theorem}\label{thm:dilation-geometry}
Let $\mu$ be a normalized OVM on $X$ admitting a cyclic dilation. In each assertion below, $Z$ is an order-continuous K\"othe function space.
\begin{enumerate}
\item $\ell_1\not\hookrightarrow X$ if and only if $\mu$ has a cyclic dilation $(Z,P,S,T)$ with $\ell_1\not\hookrightarrow Z$ and $P$ shrinking.
\item $c_0\not\hookrightarrow X$ if and only if $\mu$ has a cyclic dilation $(Z,P,S,T)$ with $Z$ a KB-space and $P$ boundedly complete.
\item $X$ is reflexive if and only if $\mu$ has a cyclic dilation $(Z,P,S,T)$ with $Z$ reflexive.
\end{enumerate}
\end{theorem}
\endgroup

\begin{proof}
{Fix any dilation $(L_1(m),Q,I_m,T_0)$ of $\mu$, where $Q(E)h=\1_Eh$, and let $Z_0=L_1(m)$. Such a dilation exists by Proposition~\ref{prop:vector-measure-dilation}. If $c_0\not\hookrightarrow X$, the proof of Corollary~\ref{cor:review-minimal-kb} shows that $Z_0$ is a KB-space. Its multiplication measure $Q$ is boundedly complete, giving the existence assertion in (2).}

Suppose now that $\ell_1\not\hookrightarrow X$. Let $K=T_0(B_X)\cup\{\1_\Omega\}$, and let $W$ be its closed absolutely convex multiplier hull as in Lemma~\ref{lem:multiplier-hull}. The inclusion $T_0(B_X)\subseteq W$ will keep the analysis operator bounded in the interpolation space. Multiplier invariance will allow the multiplication projections to act on the interpolation space, and $\1_\Omega$ will be a cyclic vector. By Rosenthal's theorem, $K$ is weakly precompact. By Lemma~\ref{lem:multiplier-hull}, $W$ is weakly precompact. If $X$ is reflexive, then $Q$ is boundedly complete by the preceding paragraph, so $W$ is weakly compact. For $n\geq1$, let $p_n$ be the Minkowski functional of $U_n=2^nW+2^{-n}B_{Z_0}$. Define
\[
F=\left\{h\in Z_0:\sum_{n=1}^{\infty}p_n(h)^2<\infty\right\}
\quad\text{with norm}\quad \norm{h}_F=\left(\sum_{n=1}^{\infty}p_n(h)^2\right)^{1/2}.
\]
By the Davis--Figiel--Johnson--Pe\l czy\'nski (DFJP) construction~\cite[Lemma~1]{DFJP1974}, $F$ is complete and its inclusion $J\colon F\to Z_0$ is bounded, and $F$ is reflexive if $W$ is weakly compact. Since $W$ is weakly precompact, $F$ contains no copy of $\ell_1$ by~\cite[Section~4]{AvilesGuiraoLajaraRodriguezTradacete2016}. The multiplier invariance of $W$ makes each $U_n$ solid, so these multipliers act contractively on $F$. The lattice DFJP construction shows that $F$ is an order-continuous K\"othe function space in $Z_0$, as in~\cite[proof of Theorem~4.3 and Remark~4.4]{AvilesGuiraoLajaraRodriguezTradacete2016}. For $h\in W$, we have $p_n(h)\leq2^{-n}$, so $\norm{h}_F^2\leq\sum_n4^{-n}=1/3$ and $W\subseteq B_F$. By order continuity, the multiplication measure $Q_F(E)h=\1_Eh$ is strongly countably additive. The truncations $h\1_{\{\abs{h}\leq k\}}$ converge to $h$ in $F$. Each truncation is a uniform limit of simple functions, and these approximations converge in $F$ because $\1_\Omega\in W\subseteq F$. Therefore $\1_\Omega$ is cyclic.

Define $T_Fx=T_0x$ and $S_F=I_mJ$. Since $T_0(B_X)\subseteq W\subseteq B_F$, the operator $T_F\colon X\to F$ is bounded. The identities $JT_F=T_0$ and $JQ_F(E)=Q(E)J$ imply
\begin{align*}
S_FT_F&=I_mJT_F=I_mT_0=I_X,
\end{align*}
and
\begin{align*}
S_FQ_F(E)T_F&=I_mQ(E)JT_F=I_mQ(E)T_0=\mu(E).
\end{align*}
Thus $(F,Q_F,S_F,T_F)$ is a cyclic dilation, and $Q_F$ is shrinking by Proposition~\ref{prop:no-l1-general}. This proves the existence assertions in (1) and (3).

Conversely, $ST=I_X$ identifies $X$ with a complemented subspace of every dilation space. By Theorem~\ref{thm:representation-james}, a shrinking or boundedly complete cyclic dilation excludes $\ell_1$ or $c_0$, respectively, from $X$. A reflexive dilation implies that $X$ is reflexive.
\end{proof}

\begingroup
The following theorem characterizes reflexive dilations and dilations on spaces containing no copy of $\ell_1$, without assuming cyclicity.

\begin{theorem}\label{prop:multiplier-orbit-criterion}
Let $\mu$ be a normalized OVM on $X$, and let $(Z_\alpha,P_\alpha,S_\alpha,T_\alpha)$ be its minimal dilation. Let $\pi_\alpha$ be the multiplier representation of $P_\alpha$ from Remark~\ref{rem:ovm-bounds}, and define
\[
K_\alpha(\mu)=\{\pi_\alpha(a)T_\alpha x:\norm{a}_\infty\leq1,\ \norm{x}\leq1\}.
\]
Then we have
\begin{enumerate}
\item $\mu$ admits a reflexive dilation if and only if $K_\alpha(\mu)$ is relatively weakly compact in $Z_\alpha$.
\item $\mu$ admits a dilation on a space containing no copy of $\ell_1$ if and only if $K_\alpha(\mu)$ is weakly precompact in $Z_\alpha$.
\end{enumerate}
\end{theorem}

\begin{proof}
Let $(Z,P,S,T)$ be a dilation, and let $W=\cspan\{P(E)Tx:E\in\Sigma,\ x\in X\}$. As in the proof of Proposition~\ref{prop:minimal-cyclic}, the formula $(R_0z)(E)=SP(E)z$ defines a bounded map $R_0\colon W\to Z_\alpha$. Let $\pi_P$ be the multiplier representation of $P$. The identities $R_0T=T_\alpha$ and $R_0\pi_P(a)|_W=\pi_\alpha(a)R_0$ show
\[
K_\alpha(\mu)=R_0\{\pi_P(a)Tx:\norm{a}_\infty\leq1,\ \norm{x}\leq1\}.
\]
If $Z$ is reflexive, the bounded multiplier orbit in $W$ is relatively weakly compact. If $\ell_1\not\hookrightarrow Z$, it is weakly precompact by Rosenthal's theorem. Since $R_0$ is bounded, $K_\alpha(\mu)$ is relatively weakly compact in the first case and weakly precompact in the second.

Conversely, suppose that $K_\alpha(\mu)$ has the compactness property in the corresponding assertion, and let $C=\overline{\operatorname{aconv}}K_\alpha(\mu)$. In (1), $C$ is weakly compact by Krein's theorem. In (2), it is weakly precompact, as recalled in~\cite[Section~4]{AvilesGuiraoLajaraRodriguezTradacete2016}. Let $p_n$ and $F$ be given by the interpolation construction in the proof of Theorem~\ref{thm:dilation-geometry}, with $Z_0=Z_\alpha$ and $W=C$. The inclusion $J\colon F\to Z_\alpha$ is bounded and $C\subseteq B_F$. The space $F$ is reflexive in (1) by~\cite[Lemma~1]{DFJP1974}, and contains no copy of $\ell_1$ in (2) by~\cite[Section~4]{AvilesGuiraoLajaraRodriguezTradacete2016}. Each $P_\alpha(E)$ is contractive on $Z_\alpha$ and preserves $C$, so it is contractive for every $p_n$ and hence on $F$. If $A_k\downarrow\varnothing$ and $z\in F$, then
\[
p_n(P_\alpha(A_k)z)\leq2^n\norm{P_\alpha(A_k)z}_\alpha\longrightarrow0
\quad\text{and}\quad p_n(P_\alpha(A_k)z)\leq p_n(z).
\]
By dominated convergence in $\sum_np_n(z)^2$, we have $\norm{P_\alpha(A_k)z}_F\to0$. Thus $P_F(E)=P_\alpha(E)|_F$ is strongly countably additive. Define $T_Fx=T_\alpha x$ and $S_F=S_\alpha J$. Since $T_\alpha(B_X)\subseteq C\subseteq B_F$, $T_F$ is bounded. By the dilation identities, we obtain$S_FT_F=I_X$ and $S_FP_F(E)T_F=\mu(E)$.
\end{proof}
\par\endgroup

\begingroup
\begin{remark}\label{rem:scalar-minimal}
The minimal $\alpha$-dilation need not be reflexive, even when $X$ is reflexive. Let $\lambda$ be Lebesgue measure on $[0,1]$, and define $\mu(E)a=\lambda(E)a$ on $\R$. For a simple real function $f$, let $\eta_f(E)=\int_Ef\,\rmd\lambda$. Then
\[
\norm{\eta_f}_\alpha=\sup_E\abs{\int_Ef\,\rmd\lambda}=\frac{1}{2}\left(\norm{f}_1+\abs{\int_0^1f\,\rmd\lambda}\right).
\]
Hence $Z_\alpha$ is isomorphic to $L_1[0,1]$ and is not reflexive, whereas $Ta=a\1$, $Sh=\int_0^1h\,\rmd\lambda$ and $P(E)h=\1_Eh$ give a reflexive cyclic dilation on $L_2[0,1]$. The measure $\mu$ is shrinking, whereas $P_\alpha$ is not by Theorem~\ref{thm:representation-james}. For the cyclic vector $u=\eta_{\1}$, we have $m=\lambda$ in Proposition~\ref{prop:vector-measure-dilation}, and taking a nonzero simple $f$ with $\int_0^1f\,\rmd\lambda=0$ shows that its lower constant $1/2$ is optimal.
\end{remark}
\endgroup

\subsection[Reflexive dilations under local variation bounds]{Reflexive dilations under local variation bounds}

The \emph{operator-norm variation} of an OVM $\mu$ is the positive measure
\[
\abs{\mu}(E)
=\sup\Big\{\sum_{j=1}^{m}\norm{\mu(E_j)}:
(E_j)_{j=1}^{m}\text{ is a measurable partition of }E\Big\},
\]
where the value $+\infty$ is allowed.

Han, Larson and Liu~\cite[Theorem~1.2]{HanLarsonLiu2018} constructed projection-valued dilations preserving bounded $p$-variation. Let $\mu_x(E)=\mu(E)x$ for $x\in X$. The quantity $\sup_{\norm{x}\leq1}\abs{\mu_x}(E)$ is the strong variation on $E$, as in~\cite[Section~2.6]{HanLarsonLiuLiu2014}. Here we obtain reflexive dilations from local finiteness of the orbit variations.

\begin{theorem}\label{thm:variation-reflexive-dilation}
Let $\mu$ be a normalized OVM on $X$. Suppose that there is a countable measurable partition $(B_j)$ of $\Omega$ such that $\abs{\mu_x}(B_j)<\infty$ for every $x\in X$ and $j\in\N$. Then $\mu$ admits a reflexive dilation if and only if $X$ is reflexive. In particular, this holds whenever $\abs{\mu}$ is $\sigma$-finite.
\end{theorem}

\begin{proof}
Every dilation contains a complemented copy of $X$, so only sufficiency requires proof. Suppose that $X$ is reflexive. For each $j$, the map $x\mapsto\abs{\mu_x}(B_j)$ is a finite lower semicontinuous seminorm, since it is the supremum of $x\mapsto\sum_i\norm{\mu(A_i)x}$ over finite measurable partitions $(A_i)$ of $B_j$. By the Baire category theorem, there is $C_j<\infty$ such that $\abs{\mu_x}(B_j)\leq C_j\norm{x}$ for every $x\in X$.
\begingroup
Let $K=K_\alpha(\mu)$ be the multiplier orbit in Theorem~\ref{prop:multiplier-orbit-criterion}. It is bounded by~\eqref{eq:multiplier-bound}. By Theorem~\ref{prop:multiplier-orbit-criterion}, it suffices to prove that $K$ is relatively weakly compact. We first prove that each $P_\alpha(B_j)K$ is relatively weakly compact, then control the scalar tails and combine the local limits using bounded completeness of $P_\alpha$.
\par\endgroup

Fix $j$, and define $V_j=\cspan\{A\mapsto\mu(A\cap E\cap B_j)x:E\in\Sigma,\ x\in X\}$, where the closure is taken in the variation norm $\norm{m}_{V_j}=\abs{m}(\Omega)$. The operator $T_j\colon X\to V_j$, defined by $(T_jx)(A)=\mu(A\cap B_j)x$, satisfies $\norm{T_j}\leq C_j$. Moreover, $\norm{am}_{V_j}\leq\norm{a}_\infty\norm{m}_{V_j}$. Since $\norm{m}_\alpha\leq\abs{m}(\Omega)$ on the elementary measures, their inclusion extends to a contraction $R_j\colon V_j\to Z_\alpha$.

We claim that $K_j=\{aT_jx:\norm{a}_\infty\leq1,\ \norm{x}\leq1\}$ is relatively weakly compact in $V_j$. Let $(a_nm_n)\subseteq K_j$, where $m_n=T_jx_n$, and define the finite positive measure
\[
\nu=\sum_{n=1}^{\infty}\frac{2^{-n}\abs{m_n}}{1+\abs{m_n}(\Omega)}.
\]
Since $X$ is reflexive, $T_j(B_X)$ is relatively weakly compact. The $\nu$-absolutely continuous vector measures form a norm-closed, hence weakly closed, subspace of the space of vector measures of bounded variation. By the Radon--Nikod\'ym property of $X$, this subspace is isometric to $L_1(\nu;X)$. Consequently, the densities $h_n=\rmd m_n/\rmd\nu$ form a relatively weakly compact subset of $L_1(\nu;X)$ and are uniformly integrable. Since $\norm{a_nh_n(t)}\leq\norm{h_n(t)}$ and $X$ is reflexive, the sequence $(a_nh_n)$ is relatively weakly compact in $L_1(\nu;X)$; see~\cite{Diestel1977}. Its image $(a_nm_n)$ therefore has a weakly convergent subsequence, whose limit belongs to the closed subspace $V_j$. By the Eberlein--\v Smulian theorem, $K_j$ is relatively weakly compact. Since $P_\alpha(B_j)\pi_\alpha(a)T_\alpha x=R_j(aT_jx)$, the set $P_\alpha(B_j)K$ is relatively weakly compact.

We next estimate the scalar tails uniformly on $K$. For $f\in Z_\alpha^*$, define $\nu_f(A)=T_\alpha^*P_\alpha(A)^*f\in X^*$. This measure is weak-star countably additive. Since $X$ is reflexive, it is weakly countably additive and therefore norm countably additive by the Orlicz--Pettis theorem. Let $D_N=\bigcup_{j>N}B_j$.
\begingroup
For $x\in B_X$, define the scalar measure $\nu_{f,x}(A)=\nu_f(A)(x)$. If $\norm{a}_\infty\leq1$, then
\[
f(P_\alpha(D_N)\pi_\alpha(a)T_\alpha x)=\int_{D_N}a\,\rmd\nu_{f,x}.
\]
The total variation of a real scalar measure satisfies
\[
\abs{\nu_{f,x}}(D_N)\leq2\sup_{A\subseteq D_N}\abs{\nu_{f,x}(A)}\leq2\sup_{A\subseteq D_N}\norm{\nu_f(A)}.
\]
By taking suprema over $a$ and $x$ and by continuity of the semivariation of $\nu_f$, we obtain
\begingroup
\[
\sup_{z\in K}\abs{f(P_\alpha(D_N)z)}
\leq2\sup_{A\subseteq D_N}\norm{\nu_f(A)}
\longrightarrow0.
\]
\endgroup
\endgroup

We now combine the local weak limits. Let $(z_k)\subseteq K$. By the Eberlein--\v Smulian theorem and a diagonal selection, pass to a subsequence such that $P_\alpha(B_j)z_k\xrightarrow{w} w_j$ for every $j$. The vectors $w_j$ belong to $P_\alpha(B_j)Z_\alpha$, and
\[
\norm{\sum_{j=1}^Nw_j}_\alpha
\leq\liminf_k\norm{P_\alpha(\Omega\setminus D_N)z_k}_\alpha
\leq\sup_k\norm{z_k}_\alpha.
\]
By Proposition~\ref{prop:minimal-bounded-complete}, the series $\sum_jw_j$ converges to some $z\in Z_\alpha$. Since $P_\alpha(B_j)z=w_j$, we have $P_\alpha(\Omega\setminus D_N)(z_k-z)=\sum_{j=1}^N(P_\alpha(B_j)z_k-w_j)\xrightarrow{w}0$ as $k\to\infty$ for each $N$. Thus $\limsup_{k\to\infty}\abs{f(z_k-z)}\leq\sup_{u\in K}\abs{f(P_\alpha(D_N)u)}+\abs{f(P_\alpha(D_N)z)}$ for every $f\in Z_\alpha^*$. By the preceding tail estimate and strong countable additivity of $P_\alpha$, both terms on the right tend to zero as $N\to\infty$, so $z_k\xrightarrow{w}z$. Thus $K$ is relatively weakly compact.

Theorem~\ref{prop:multiplier-orbit-criterion} now gives a reflexive dilation.
\end{proof}

The next example shows that bounded $2$-variation does not ensure a reflexive dilation.

\begin{example}\label{ex:no-reflexive-dilation}
There is a normalized, operator-norm countably additive measure on $\ell_2$ with bounded $2$-variation and no reflexive dilation. Let $(r_n)$ be the independent Rademacher functions on $[0,1]$, with Lebesgue measure $\lambda$, and define
\[
(\mu(A)x)_n=\sum_{N\geq n}x_N\int_A r_nr_N\,\rmd\lambda.
\]
For finitely supported $x,f\in\ell_2$, the scalar measure $f\mu(\cdot)x$ has density
\[
g_{x,f}=\sum_nx_nf_n+\sum_{n<N}x_Nf_nr_nr_N.
\]
By Walsh orthogonality,
\begin{align*}
\norm{g_{x,f}}_2^2
&=\abs{\sum_nx_nf_n}^2+\sum_{n<N}\abs{x_Nf_n}^2
  \leq2\norm{x}_2^2\norm{f}_2^2,
\end{align*}
and
\begin{align*}
\abs{f(\mu(A)x)}
&\leq\lambda(A)^{1/2}\norm{g_{x,f}}_2
  \leq\sqrt2\lambda(A)^{1/2}\norm{x}_2\norm{f}_2.
\end{align*}
Thus the formula extends to $\mu(A)\in\cB(\ell_2)$, with $\norm{\mu(A)}\leq\sqrt{2\lambda(A)}$. It follows that $\mu$ is operator-norm countably additive, and orthogonality implies $\mu([0,1])=I_{\ell_2}$. For every finite measurable partition $(A_i)$ of $[0,1]$, we also have $\sum_i\norm{\mu(A_i)x}_2^2\leq2\sum_i\lambda(A_i)\norm{x}_2^2=2\norm{x}_2^2$, so $\mu$ has bounded $2$-variation.

Suppose that $(F,P,S,T)$ is a reflexive dilation, and let $\pi$ be the multiplier representation of $P$. The sequences $u_N=\pi(r_N)Te_N$ and $v_j=\pi(r_j)^*S^*e_j^*$ are bounded in $F$ and $F^*$, respectively, and
\[
v_j(u_N)=e_j^*\left(\int r_jr_N\,\rmd\mu\,e_N\right)
=\begin{cases}1,&j\leq N,\\0,&j>N.\end{cases}
\]
By reflexivity, pass to subsequences with $u_{N_k}\xrightarrow{w} u$ and $v_{j_l}\xrightarrow{w} v$. For every $j$, we have $v_j(u)=1$, whereas $v(u_N)=0$ for every $N$. Thus $v(u)=1$ and $v(u)=0$, a contradiction.

Finally, $\mu(\cdot)e_N$ has Bochner density $r_N\sum_{n=1}^Nr_ne_n$, whose norm is $\sqrt N$. Consequently, $\abs{\mu}(A)\geq\abs{\mu(\cdot)e_N}(A)=\sqrt N\lambda(A)$ for every $N$. If $\lambda(A)>0$, the orbit variations are unbounded on $B_{\ell_2}$. By the Baire category argument in Theorem~\ref{thm:variation-reflexive-dilation}, they cannot all be finite on $A$. Thus no countable measurable partition satisfies its hypothesis, and $\abs{\mu}$ is not $\sigma$-finite.
\end{example}

\section{Schauder frames}\label{sec:frames}

We first ask when a given OVM admits a continuous Schauder frame representation on its measurable space. We then study the geometry of associated spaces for a fixed frame, retaining its original analysis coefficients.

Throughout this section, $X$ is separable.

\subsection{From frames to cyclic dilations}

Let $(\Omega,\Sigma,\eta)$ be a $\sigma$-finite measure space.

\begin{definition}
A weakly measurable map $H\colon\Omega\to X$ is \emph{Pettis integrable} if $g\circ H\in L_1(\eta)$ for every $g\in X^*$, and for each $E\in\Sigma$ there is a vector $x_E\in X$ such that
\[
g(x_E)=\int_Eg(H(t))\,\rmd\eta(t)\quad(g\in X^*).
\]
We denote the vector $x_E$ by $(P)\!\int_EH\,\rmd\eta$. A family $(x_t,f_t)_{t\in\Omega}\subseteq X\times X^*$ is a \emph{continuous Schauder frame} if $t\mapsto x_t$ is weakly measurable, $t\mapsto f_t$ is weak-star measurable, and
\[
x=(P)\!\int_\Omega f_t(x)x_t\,\rmd\eta(t)\quad(x\in X).
\]
\end{definition}

This is the convention of~\cite[Definitions~2.2 and~2.3]{EisnerFreeman2020}. For counting measure on $\N$, a continuous Schauder frame is an unconditional Schauder frame.

Li, Li and Han~\cite[Theorem~4.3]{LiLiHan2016} characterize continuous framings through a $g$-K\"othe realization with local representability conditions on the adjoint synthesis operator. Their Proposition~4.8 gives compact reconstruction operators on a countable measurable partition when both frame components are strongly measurable. Under the present separability assumption, the following proposition only requires weak-star measurability of $t\mapsto f_t$. It also gives $\sigma$-finite operator-norm variation and a cyclic vector-measure dilation. The following proposition gives the necessary conditions for the frame representation in Theorem~\ref{thm:continuous-frame-representation}.

\begingroup
\begin{proposition}\label{prop:frame-dilation}
A continuous Schauder frame induces a normalized OVM
\[
K(E)x=(P)\!\int_E f_t(x)x_t\,\rmd\eta(t).
\]
It admits a cyclic dilation on $L_1(m)$ for some countably additive vector measure $m\colon\Sigma\to X$. There is a countable measurable partition $(D_n)$ of $\Omega$ with $\abs{K}(D_n)<\infty$. For this partition, $K(E)$ is compact whenever $E\subseteq D_n$ is measurable.
\end{proposition}
\endgroup

\begin{proof}
For $g\in X^*$, Pettis integrability defines $T_g:X\to L_1(\eta)$ by $T_gx(t)=f_t(x)g(x_t)$. This operator has closed graph: if $x_j\to x$ and $T_gx_j\to h$ in $L_1$, an a.e.\ convergent subsequence identifies $h=T_gx$. Hence $T_g$ is bounded, and $\abs{g(K(E)x)}\leq\norm{T_g}\norm{x}$. By the uniform boundedness principle, $K(E)\in\cB(X)$. Scalar integration, reconstruction and Remark~\ref{rem:ovm-bounds} show that $K$ is a normalized OVM.

By separability, $t\mapsto x_t$ is strongly measurable and $t\mapsto\norm{f_t}$ is measurable. Choose a measurable partition $(D_n)$ such that $\eta(D_n)<\infty$ and both norms are bounded on each $D_n$. By the Pettis integral estimate, we have
\[
\abs{K}(E)\leq\int_E\norm{x_t}\norm{f_t}\,\rmd\eta(t),
\]
so $\abs{K}(D_n)<\infty$. On each $D_n$, choose finite-valued simple functions $s_{n,j}$ with
\[
\int_{D_n}\norm{x_t-s_{n,j}(t)}\norm{f_t}\,\rmd\eta(t)\to0.
\]
For $E\subseteq D_n$, the finite-rank operators $K_{n,j}(E)x=\int_E f_t(x)s_{n,j}(t)\,\rmd\eta(t)$ satisfy
\[
\norm{K(E)-K_{n,j}(E)}\leq\int_{D_n}\norm{x_t-s_{n,j}(t)}\norm{f_t}\,\rmd\eta(t)\longrightarrow0.
\]
Thus $K(E)$ is compact.

On $D_n$, define
\[
w(t)=2^{-n}/[(1+\eta(D_n))(1+\norm{x_t})].
\]
Then $w>0$ and $\int_\Omega w\norm{x_t}\,\rmd\eta\leq1$, so $m(E)=\int_Ew(t)x_t\,\rmd\eta(t)$ defines a countably additive vector measure by Bochner integration. Define $(T_0x)(t)=f_t(x)/w(t)$ and $Q(E)h=\1_Eh$. For every $x\in X$ and $g\in X^*$,
\begin{align*}
\int_\Omega\abs{T_0x}\,\rmd\abs{gm}
&=\int_\Omega\abs{f_t(x)g(x_t)}\,\rmd\eta(t)<\infty
  \quad\text{and}\quad \int_E T_0x\,\rmd(gm)=g(K(E)x).
\end{align*}
Thus $T_0x\in L_1(m)$ and $I_mQ(E)T_0x=K(E)x$. By scalar total variation,
\begin{align*}
\norm{T_0x}_{L_1(m)}
&=\sup_{g\in B_{X^*}}\abs{gK(\cdot)x}(\Omega)\\
&\leq c\sup_{g\in B_{X^*}}\sup_{E\in\Sigma}\abs{g(K(E)x)}
\leq cM_K\norm{x},
\end{align*}
where $c=2$ over $\R$ and $c=4$ over $\mathbb C$. Finally, $I_mT_0=I_X$, and $\1_\Omega$ is cyclic for $Q$.
\end{proof}

The simple-function approximation in the preceding proof also gives an unconditional Schauder frame by averaging the coefficient functionals over a measurable partition and applying an invertible correction to the synthesis vectors.

\begin{corollary}\label{cor:frame-discretization}
Every separable Banach space admitting a continuous Schauder frame admits an unconditional Schauder frame.
\end{corollary}

\begin{proof}
Use the partition and simple functions from the preceding proof. For $0<\varepsilon<1$, choose one $s_n=s_{n,j(n)}$ on each $D_n$ with weighted error less than $2^{-n}\varepsilon$. Enumerate their cells as $(B_j)$, with corresponding values $y_j$, and define $g_j(x)=\int_{B_j}f_t(x)\,\rmd\eta(t)$. Since $\norm{g_j}\leq\int_{B_j}\norm{f_t}\,\rmd\eta<\infty$, we have $g_j\in X^*$. Moreover,
\begin{align*}
\sum_j\norm{y_j\otimes g_j-K(B_j)}
&\leq\sum_j\int_{B_j}\norm{y_j-x_t}\norm{f_t}\,\rmd\eta(t)
  <\varepsilon.
\end{align*}
Since $\sum_jK(B_j)x=x$ unconditionally, the series $Ax=\sum_jg_j(x)y_j$ also converges unconditionally, and $\norm{A-I_X}<\varepsilon$. Thus $A$ is invertible and $(A^{-1}y_j,g_j)$ is an unconditional Schauder frame.
\end{proof}

\subsection{Frame representations of OVMs}

\begingroup
An \emph{unconditional Schauder decomposition} of a Banach space $Z$ is a family $(Z_n)$ of closed subspaces such that every $z\in Z$ has a unique unconditionally convergent expansion $z=\sum_nz_n$, with $z_n\in Z_n$. It is \emph{finite-dimensional} if every $Z_n$ is finite-dimensional.

\begin{proposition}\label{prop:discrete-frame-representation}
Let $\mu\colon2^{\N}\to\cB(X)$ be a normalized OVM. Then $X$ is reflexive if and only if $\mu$ has a reflexive dilation $(Z,P,S,T)$. In this case, the dilation may be chosen so that the nonzero spaces $Z_n=P(\{n\})Z$ form an unconditional Schauder decomposition and $\dim Z_n=\operatorname{rank}\mu(\{n\})$ whenever $\mu(\{n\})$ has finite rank.

For each $r\in\N$, the following are equivalent without assuming that $X$ is reflexive:
\begin{enumerate}
\item $\operatorname{rank}\mu(\{n\})\leq r$ for every $n\in\N$.
\item $\mu$ has a dilation $(Z,P,S,T)$ with $ Z=\cspan\{P(E)u_j:E\subseteq\N,\ 1\leq j\leq r\} $ for some $u_1,\ldots,u_r\in Z$.
\item the minimal $\alpha$-dilation is generated by at most $r$ vectors in the same sense.
\end{enumerate}
For $r=1$, these conditions are equivalent to the existence of an unconditional Schauder frame $(x_n,f_n)_{n\in\N}$ inducing $\mu$.
\end{proposition}

\begin{proof}
For $z\in P_\alpha(\{n\})Z_\alpha$, we have
\[
z(E)=\1_E(n)z(\{n\})\quad\text{and}\quad\norm{z}_\alpha=\norm{z(\{n\})}.
\]
Evaluation at $\{n\}$ is therefore an isometry onto $\overline{\mu(\{n\})X}$, since the projected elementary measures are dense and their evaluations range over $\mu(\{n\})X$.

Since $\abs{\mu_x}(\{n\})=\norm{\mu(\{n\})x}<\infty$, the reflexive dilation assertion follows from Theorem~\ref{thm:variation-reflexive-dilation}. Its construction gives a subspace $F\subseteq Z_\alpha$, and we define $Z_n=P_\alpha(\{n\})F$. By strong countable additivity, the nonzero spaces $Z_n$ form an unconditional Schauder decomposition of $F$. Each $Z_n$ contains $P_\alpha(\{n\})T_\alpha X$, so the preceding identification gives the asserted dimension when the rank is finite.

Under (2), $P(\{n\})Z=\spn\{P(\{n\})u_j:1\leq j\leq r\}$, so $\mu(\{n\})=SP(\{n\})T$ has rank at most $r$. Conversely, suppose (1), and let $d_n=\operatorname{rank}\mu(\{n\})$. Choose a normalized basis $v_{n,1},\ldots,v_{n,d_n}$ of $P_\alpha(\{n\})Z_\alpha$, and define $v_{n,j}=0$ for $d_n<j\leq r$. The absolutely convergent sums $u_j=\sum_n2^{-n}v_{n,j}$ satisfy $P_\alpha(\{n\})u_j=2^{-n}v_{n,j}$. Hence their closed invariant span contains every coordinate space and equals $Z_\alpha$ by strong countable additivity. This proves (3), which implies (2).

When $r=1$, factor $\mu(\{n\})=x_n\otimes f_n$, taking both factors zero if the operator vanishes. Strong countable additivity and normalization give the required unconditional frame expansion. Conversely, every such frame has a cyclic dilation by Proposition~\ref{prop:frame-dilation}.
\end{proof}
\endgroup

Li and Li~\cite[Theorem~2.6]{LiLi2017} obtain operator-valued densities under the Radon--Nikod\'ym assumption. Under the hypotheses below, we first derive this property from cyclicity and $\sigma$-finite operator-norm variation.

\begin{proposition}\label{prop:variation-rnp}
Let $X$ be a separable Banach space containing no copy of $c_0$. If a normalized OVM $\mu$ on $X$ admits a cyclic dilation and has $\sigma$-finite operator-norm variation, then $X$ has the Radon--Nikod\'ym property.
\end{proposition}

\begin{proof}
By Proposition~\ref{prop:vector-measure-dilation}, we have $\mu(E)=I_mQ(E)T_0$ and $I_mT_0=I_X$. {Since $c_0\not\hookrightarrow X$, the space $L_1(m)$ is KB by the proof of Corollary~\ref{cor:review-minimal-kb}.} We first show that $\mu(D)$ is weakly compact whenever $\abs{\mu}(D)<\infty$. Choose a dense sequence $(d_j)$ in $B_X$, and let $h_j$ be a measurable representative of $T_0d_j$. Let $H_k=\1_D\max_{1\leq j\leq k}\abs{h_j}$. For fixed $k$, partition $D$ into measurable sets $E_1,\ldots,E_k$, assigning each point to the first index attaining the maximum. Since $d_j\in B_X$, scalar variation satisfies $\abs{f\mu_{d_j}}(B)\leq \abs{\mu}(B)$ for every $f\in B_{X^*}$ and measurable $B$. Since $f\mu_{d_j}=h_j\cdot(fm)$, for $f\in B_{X^*}$ we have
\begin{align*}
\int_\Omega H_k\,\rmd\abs{fm}
&=\sum_{j=1}^{k}\int_{E_j}\abs{h_j}\,\rmd\abs{fm}\\*
&=\sum_{j=1}^{k}\abs{f\mu_{d_j}}(E_j)
\leq\sum_{j=1}^{k}\abs{\mu}(E_j)=\abs{\mu}(D).
\end{align*}
Thus $\norm{H_k}_{L_1(m)}\leq \abs{\mu}(D)$. By the KB property, there is $H\in L_1(m)$ such that $H_k\uparrow H$ in norm. For each $j$, the inequality $\1_D\abs{h_j}\leq H_k$ holds whenever $k\geq j$, so $Q(D)T_0d_j\in[-H,H]$. Since $Q(D)T_0$ is bounded and this order interval is norm closed, we have by density $Q(D)T_0(B_X)\subseteq[-H,H]$. The norm of $L_1(m)$ is order continuous, so $[-H,H]$ is weakly compact. Hence $\mu(D)=I_mQ(D)T_0$ is weakly compact.

Choose a measurable partition $(D_n)$ of $\Omega$ with $\abs{\mu}(D_n)<\infty$, and let $A_n=\mu(D_n)$. By strong countable additivity and normalization, we have $x=\sum_nA_nx$ unconditionally. By Remark~\ref{rem:ovm-bounds},
\[
\Big\|\sum_{n\in J}\varepsilon_nA_n\Big\|\leq 2M_\mu
\quad(J\subseteq\N\text{ finite},\ \varepsilon_n\in\{-1,1\}).
\]

Let $\nu\ll\lambda$ be an $X$-valued measure of bounded variation on $(\Lambda,\mathcal A)$, where $\lambda$ is a finite positive measure, and let $r=\rmd\abs{\nu}/\rmd\lambda$. By DFJP factorization~\cite{DFJP1974}, each $A_n$ factors through a reflexive space. Since reflexive spaces have the Radon--Nikod\'ym property, $A_n\nu$ has a Bochner density $g_n$. For every finite $J$ and choice of signs, we have $\big|(\sum_{n\in J}\varepsilon_nA_n)\nu\big|\leq 2M_\mu\abs{\nu}$. The measure on the left has density $\sum_{n\in J}\varepsilon_ng_n$. By the variation formula for Bochner integrals, we have $\big\|\sum_{n\in J}\varepsilon_ng_n(t)\big\|\leq 2M_\mu r(t)$ a.e. There are only countably many finite sets and sign choices, so these inequalities hold outside one common null set.

At each remaining point, the bounded finite signed sums imply that $\sum_ng_n(t)$ is weakly unconditionally Cauchy. By the Bessaga--Pe\l czy\'nski theorem~\cite{BessagaPelczynski1958}, the series converges unconditionally in norm, since $X$ contains no copy of $c_0$. Denote the sum by $g(t)$, and define $g=0$ on the exceptional set. The partial sums and their limit $g$ are strongly measurable and bounded in norm by $2M_\mu r$. Since $r\in L_1(\lambda)$, we have by dominated convergence
\[
\int_Bg\,\rmd\lambda
=\lim_{N\to\infty}\sum_{n=1}^{N}A_n\nu(B)=\nu(B)
\quad(B\in\mathcal A).
\]
Thus $X$ has the Radon--Nikod\'ym property.
\end{proof}

The Radon--Nikod\'ym property obtained above provides an operator-valued density. Cyclicity then forces this density to have rank at most one, yielding a continuous Schauder frame representation of the given OVM on the original measurable space.

\Needspace{6\baselineskip}
\begin{theorem}\label{thm:continuous-frame-representation}
Let $X$ be a separable Banach space containing no copy of $c_0$, and let $\mu\colon\Sigma\to\cB(X)$ be a normalized OVM. The following are equivalent:
\begin{enumerate}
\item The measure $\mu$ is induced by a continuous Schauder frame $(x_t,f_t)_{t\in\Omega}$ over $(\Omega,\Sigma,\lambda)$ for some $\sigma$-finite positive measure $\lambda$.
\item The measure $\mu$ admits a cyclic dilation and its operator-norm variation $\abs{\mu}$ is $\sigma$-finite.
\end{enumerate}
In (1), one may take $\lambda=\abs{\mu}$ and $\norm{x_t}\norm{f_t}=1$ for $\abs{\mu}$-almost every $t$. If $X^*$ is separable, both maps are strongly measurable.
\end{theorem}

\begin{proof}
By Proposition~\ref{prop:frame-dilation}, (1) implies (2). Conversely, suppose that (2) holds. We first obtain an operator-valued density with norm at most one a.e. Proposition~\ref{prop:variation-rnp} gives the Radon--Nikod\'ym property of $X$, so~\cite[Theorem~2.6]{LiLi2017} supplies a density $A:\Omega\to\cB(X)$ with strongly measurable, Pettis integrable orbits satisfying
\begin{equation}\label{eq:global-scalar-density}
g(\mu(E)x)=\int_Eg(A(t)x)\,\rmd\abs{\mu}(t)\quad(E\in\Sigma,\ x\in X,\ g\in X^*).
\end{equation}

Let $D\subseteq X$ be a countable dense rational linear subspace, and let $G=\{g_j:j\in\N\}\subseteq B_{X^*}$ be a countable norming family. Uniqueness of scalar Radon--Nikod\'ym derivatives and the bound $\abs{g_j\mu_d}\leq\norm{d}\abs{\mu}$ give $\abs{g_j(A(t)d)}\leq\norm{d}$ outside one common null set for all $d\in D$ and $j\in\N$. The norming property of $G$ and density of $D$ imply $\norm{A(t)}\leq1$ after defining $A=0$ on this null set. Each orbit is therefore Bochner integrable on sets of finite $\abs{\mu}$-measure, and~\eqref{eq:global-scalar-density} gives
\begin{equation}\label{eq:local-operator-density}
\mu(E)x=\int_EA(t)x\,\rmd\abs{\mu}(t)\quad(\abs{\mu}(E)<\infty,\ x\in X).
\end{equation}

We now use cyclicity to prove that $\operatorname{rank}A(t)\leq1$ a.e. By Proposition~\ref{prop:vector-measure-dilation}, we have $\mu(E)x=\int_ET_0x\,\rmd m$. For $d\in D$, choose a measurable representative $h_d$ of $T_0d$.
\begingroup
Let $\rho$ be a finite control measure for $m$. If $\rho(E)=0$, then $\mu(B)x=\int_BT_0x\,\rmd m=0$ for every measurable $B\subseteq E$ and $x\in X$. Thus $\abs{\mu}\ll\rho$, and we may define $a=\rmd\abs{\mu}/\rmd\rho$ and $k_j=\rmd(g_jm)/\rmd\rho$. Since $\abs{\mu}$ is $\sigma$-finite, there are measurable sets $\Omega_n$ covering $\Omega$ with $ \int_{\Omega_n}a\,\rmd\rho=\abs{\mu}(\Omega_n)<\infty$. Hence $a<\infty$ $\rho$-a.e. Since $\abs{\mu}\ll\rho$ and $\abs{\mu}(\{a=0\})=0$, we have $0<a<\infty$, $\abs{\mu}$-a.e.
\endgroup

By equation~\eqref{eq:global-scalar-density} and the identity $g_j\mu_d=h_d\cdot(g_jm)$, we have two scalar densities of the same measure. By uniqueness,
\[
a(t)g_j(A(t)d)=h_d(t)k_j(t)\quad(d\in D,\ j\in\N)
\]
outside one $\rho$-null set. We show that all $2\times2$ minors of $(g_j(A(t)d))_{j\in\N,\,d\in D}$ vanish $\abs{\mu}$-a.e. Multiply the identities for $(d,j)$ and $(e,k)$, and then those for $(e,j)$ and $(d,k)$. Since both products on the right equal $h_d(t)h_e(t)k_j(t)k_k(t)$, it follows that
\[
a(t)^2\big[g_j(A(t)d)g_k(A(t)e)-g_j(A(t)e)g_k(A(t)d)\big]=0.
\]
Since $0<a(t)<\infty$ outside one $\abs{\mu}$-null set, we get
\[
g_j(A(t)d)g_k(A(t)e)-g_j(A(t)e)g_k(A(t)d)=0
\]
there, simultaneously for $d,e\in D$ and $j,k\in\N$. Fix a point outside this set. If $A(t)\neq0$, then by density of $D$ and the norming property of $G$, there are $d\in D$ and $j$ with $g_j(A(t)d)\neq0$. For every $e\in D$, the preceding identities imply
\[
g_k\left(A(t)e-\frac{g_j(A(t)e)}{g_j(A(t)d)}A(t)d\right)=0
\quad(k\in\N).
\]
Since $G$ separates points, $A(t)e$ belongs to $\spn\{A(t)d\}$. By continuity, $A(t)X\subseteq\spn\{A(t)d\}$, so $\operatorname{rank}A(t)\leq1$. Define $A=0$ on the additional null set.

To choose the rank-one factors measurably, enumerate $D=\{d_i:i\in\N\}$ and fix an enumeration of $\N^2$. The functions $g_j(A(\cdot)d_i)$ are measurable, so we can partition $\{A\neq0\}$ into measurable sets according to the first pair $(i,j)$ for which $g_j(A(t)d_i)\neq0$. On the corresponding set, define
\[
x_t=A(t)d_i\quad\text{and}\quad f_t=\frac{A(t)^*g_j}{g_j(A(t)d_i)}.
\]
Define both factors to be zero where $A=0$. By the rank-one property, we have $A(t)x=f_t(x)x_t$. On each member of this countable partition, the formulas show that $x_\cdot$ is strongly measurable and $f_\cdot$ is weak-star measurable. By~\eqref{eq:global-scalar-density}, these factors represent $\mu$, and $\mu(\Omega)=I_X$ proves reconstruction.

It remains to prove $\norm{x_t}\norm{f_t}=1$ for $\abs{\mu}$-a.e. $t$. Since $X$ is separable, $t\mapsto\norm{A(t)}$ is measurable. For a finite measurable partition $(E_j)$ of a set $E$ with $\abs{\mu}(E)<\infty$, we have by~\eqref{eq:local-operator-density} that $ \sum_j\norm{\mu(E_j)} \leq\int_E\norm{A(t)}\,\rmd\abs{\mu}(t)\leq \abs{\mu}(E) $. By taking the supremum over partitions, we have equality throughout. Since $\abs{\mu}$ is $\sigma$-finite and $\norm{A(t)}\leq1$, it follows that $\norm{A(t)}=1$ a.e. Since $A(t)=x_t\otimes f_t$, we have $\norm{x_t}\norm{f_t}=1$ a.e. Finally, if $X^*$ is separable, then $t\mapsto\norm{f_t-g}$ is measurable for each $g\in X^*$, by taking the supremum over a countable dense subset of $B_X$. By approximation from a countable norm-dense subset of $X^*$, the map $f_\cdot$ is strongly measurable.
\end{proof}

\Needspace{5\baselineskip}
The following example shows that the orbit-variation condition in Theorem~\ref{thm:variation-reflexive-dilation} does not ensure a continuous frame representation.

\begin{example}\label{ex:orbit-variation-no-frame}
There is a normalized OVM on $\ell_2$ satisfying the hypothesis of Theorem~\ref{thm:variation-reflexive-dilation} which admits a cyclic reflexive dilation but no continuous Schauder frame representation on its measurable space. Let $(r_n)$ be the independent Rademacher functions on $[0,1]$, with Lebesgue measure $\lambda$, and define the isometry $R\colon\ell_2\to L_2[0,1]$ by $Rx=\sum_nx_nr_n$. On $\Omega=[0,1]\sqcup\N$, define
\[
\mu(A\sqcup B)x=e_1\int_A Rx\,\rmd\lambda+\sum_{n\in B}x_ne_n
\quad(A\subseteq[0,1],\ B\subseteq\N).
\]
Here $A$ is Lebesgue measurable and $B$ is arbitrary. Since $\int_0^1Rx\,\rmd\lambda=0$, we have $\mu(\Omega)=I_{\ell_2}$. The measure is strongly countably additive. Moreover,
\begin{align*}
\abs{\mu_x}([0,1])
&=\int_0^1\abs{Rx(t)}\,\rmd\lambda(t)
  \leq\norm{Rx}_2=\norm{x}_2.
\end{align*}
Since $\abs{\mu_x}(\{n\})=\abs{x_n}$, the partition $[0,1],\{1\},\{2\},\ldots$ satisfies the orbit-variation hypothesis.

A cyclic reflexive dilation is defined on $Z=L_2[0,1]\oplus_2\ell_2$ by
\begin{align*}
Tx&=(Rx,x),\quad S(h,z)=e_1\int_0^1h\,\rmd\lambda+z,
  \quad\text{and}\quad P(A\sqcup B)(h,z)=(\1_Ah,\1_Bz).
\end{align*}
Both $T$ and $S$ have norm at most $\sqrt2$. Moreover, $ST=I_{\ell_2}$ and $SP(E)T=\mu(E)$. The vector $(\1,(2^{-n})_n)$ is cyclic for $P$.

Let $Q_N$ be the projection onto the first $N$ coordinates of $\ell_2$. On $[0,1]$, the measure $\mu(\cdot)Q_N$ has operator-valued Bochner density $e_1\otimes\sum_{n=1}^Nr_n(t)e_n^*$, whose norm is $\sqrt N$. Hence, for every measurable $A\subseteq[0,1]$,
\[
\abs{\mu}(A)\geq\abs{\mu(\cdot)Q_N}(A)
=\sqrt N\lambda(A)\quad(N\in\N).
\]
Thus $\abs{\mu}(A)=\infty$ whenever $\lambda(A)>0$, and $\abs{\mu}$ is not $\sigma$-finite. By Proposition~\ref{prop:frame-dilation}, no continuous Schauder frame represents $\mu$ on $\Omega$. In this example, each $\mu(A)$, $A\subseteq[0,1]$, has rank at most one, but the functional $x\mapsto\sum_nx_nr_n(t)$, initially defined on $c_{00}$, has no bounded extension to $\ell_2$, since $(r_n(t))_n\notin\ell_2$.
\end{example}

\begingroup
The next example shows that the exclusion of $c_0$ in Theorem~\ref{thm:continuous-frame-representation} cannot be omitted.

\begin{example}\label{ex:c0-variation-no-frame}
Let $(r_n)$ be $\{-1,1\}$-valued Rademacher functions on $[0,1]$, with Lebesgue measure $\lambda$. Define $Vh=(\int_0^1h r_n\,\rmd\lambda)_n$ for $h\in L_2[0,1]$, and let $m(A)=V\1_A$. Bessel's inequality gives a bounded map $V\colon L_2[0,1]\to c_0$, with $V\1=0$. Since $\norm{m(A)}_\infty\leq\lambda(A)$, $m$ is norm countably additive and $\abs{m}\leq\lambda$. Its first coordinate has variation $\lambda$, so $\abs{m}=\lambda$.

On $\Omega=[0,1]\sqcup\N$, define
\[
\mu(A\sqcup B)x=x_1m(A)+\1_Bx\quad(x\in c_0).
\]
This is a normalized OVM with $\abs{\mu}=\lambda\sqcup\#$, where $\#$ denotes counting measure. On $Z=L_2[0,1]\oplus_2c_0$, define
\[
Tx=(x_1\1,x),\quad S(h,z)=Vh+z,\text{ and}\quad P(A\sqcup B)(h,z)=(\1_Ah,\1_Bz).
\]
These maps give a dilation, and $(\1,(2^{-n})_n)$ is cyclic for $P$.

Suppose that a continuous Schauder frame represents $\mu$ with respect to a $\sigma$-finite positive measure $\eta$ on $\Omega$. On $[0,1]$, let $h(t)=f_t(e_1)x_t\in c_0$, and denote its coordinates by $h_n$. The representation shows
\[
\int_Ah_n\,\rmd\eta=\int_A r_n\,\rmd\lambda\quad(n\in\N).
\]
The first coordinate yields $\lambda=\abs{h_1}\eta$ on $[0,1]$. Let $a=\abs{h_1}$. Uniqueness of scalar Radon--Nikod\'ym derivatives gives $h_n=a r_n$ simultaneously for all $n$ outside one $\eta$-null set. Thus $\abs{h_n}=a$ there. Since $h(t)\in c_0$, we have $a=0$. This contradicts $\int_{[0,1]}a\,\rmd\eta=1$.
\end{example}
\par\endgroup

\subsection{Associated spaces of a fixed frame}

Associated spaces with prescribed geometry have been studied for discrete Schauder frames. Beanland, Freeman and Liu~\cite[Theorem~1.5]{BeanlandFreemanLiu2015} prove that a Schauder frame has a shrinking associated basis if and only if it is shrinking, and has a reflexive associated space if and only if it is shrinking and boundedly complete. Beanland and Freeman~\cite{BeanlandFreeman2024} give an explicit construction of a shrinking associated basis for a given shrinking Schauder frame. For continuous Schauder frames, Eisner and Freeman~\cite[Theorem~3.11]{EisnerFreeman2020} prove that a frame of a separable Banach space is shrinking and boundedly complete if and only if the space is reflexive. Here we construct associated K\"othe function spaces for continuous Schauder frames, with multiplication projections and the geometry specified in Corollary~\ref{cor:associated-spaces}.

\begingroup
\begin{definition}\label{def:associated-space}
Let $(x_t,f_t)_{t\in\Omega}$ be a continuous Schauder frame of $X$ over $(\Omega,\Sigma,\eta)$. We say that a K\"othe function space $\mathcal E$ on $(\Omega,\Sigma,\eta)$ is an \emph{associated space} of the frame if the following hold.
\begin{enumerate}
\item The analysis operator $T\colon X\to\mathcal E$ is an isomorphic embedding, where
\[
(Tx)(t)=f_t(x)\quad(x\in X,\ t\in\Omega).
\]
\item The synthesis operator $S\colon\mathcal E\to X$ is a bounded quotient map, where
\[
Sh=(P)\!\int_\Omega h(t)x_t\,\rmd\eta(t)\quad(h\in\mathcal E).
\]
\end{enumerate}
\end{definition}
\endgroup
The frame identity gives $ST=I_X$.
\begin{corollary}\label{cor:associated-spaces}
Let $(x_t,f_t)$ be a continuous Schauder frame of $X$. Then:
\begin{enumerate}
\item $\ell_1\not\hookrightarrow X$ if and only if the frame has an order-continuous associated K\"othe function space containing no copy of $\ell_1$;
\item $c_0\not\hookrightarrow X$ if and only if the frame has a KB associated K\"othe function space;
\item $X$ is reflexive if and only if the frame has a reflexive associated K\"othe function space.
\end{enumerate}
\end{corollary}

\begin{proof}
\begingroup
Since $ST=I_X$, every associated space contains a complemented copy of $X$, so the reverse implications hold. Let $m,w,T_0$ be as in the proof of Proposition~\ref{prop:frame-dilation}, and define $A_0=\{t\in\Omega:x_t=0\}$ and $A=\Omega\setminus A_0$.

Let $(g_j)\subseteq B_{X^*}$ be a countable norming family. If $E\subseteq A$ is measurable and $m$-null, then
\[
\int_Ew(t)\abs{g_j(x_t)}\,\rmd\eta(t)=\abs{g_jm}(E)=0\quad(j\in\N).
\]
Since $w>0$, we have $g_j(x_t)=0$ for every $j$ outside one $\eta$-null subset of $E$. By the norming property, $x_t=0$ there, so $\eta(E)=0$ because $x_t\neq0$ on $A$. Conversely, every $\eta$-null set is $m$-null by the definition of $m$. Thus the $m$-null and $\eta$-null sets on $A$ coincide.

By the construction in the proof of Theorem~\ref{thm:dilation-geometry}, we obtain an order-continuous K\"othe function space $F\subseteq L_1(m)$ with $\1_A\in F$, bounded maps $T_0\colon X\to F$ and $I_m|_F\colon F\to X$, and the geometry in the corresponding assertion. In (2), we take $F=L_1(m)$. We regard $F$ as a K\"othe function space on $A$, and define $\mathcal E_A=\{wh|_A:h\in F\}$ with $\norm{k}_{\mathcal E_A}=\norm{k/w}_F$. Multiplication by $w|_A$ is an isometric lattice isomorphism from $F$ onto $\mathcal E_A$. Since $(T_0x)(t)=f_t(x)/w(t)$, the analysis and synthesis maps on $A$ are
\[
(Tx)(t)=f_t(x)\quad\text{and}\quad Sk=I_m(k/w)=(P)\!\int_Ak(t)x_t\,\rmd\eta(t).
\]
Although $x_t=0$ on $A_0$, the functionals $f_t$ may be nonzero there. We therefore adjoin a space on $A_0$ to retain the original analysis coefficients. Choose a positive measurable $a$ on $A_0$ such that $\int_{A_0}a(1+\norm{f_t}^2)\,\rmd\eta<\infty$, and define
\[
\mathcal E=\mathcal E_A\oplus_2 L_2(A_0,a\,\rmd\eta).
\]
Since
\[
\int_{A_0}\abs{f_t(x)}^2a(t)\,\rmd\eta(t)\leq\norm{x}^2\int_{A_0}a(t)\norm{f_t}^2\,\rmd\eta(t),
\]
the analysis map $x\mapsto f_\cdot(x)|_{A_0}$ is bounded into $L_2(A_0,a\,\rmd\eta)$. Since $x_t=0$ on $A_0$, the synthesis map vanishes on that summand. Thus $\mathcal E$ has the geometry in the corresponding assertion, and the full analysis and synthesis maps are bounded and satisfy $ST=I_X$ and $S\1_ET=K(E)$. By order continuity of $\mathcal E$, the multiplication measure is strongly countably additive.

Let $v=w\1_A+\1_{A_0}$. Since $\1_A$ is cyclic in $F$, we have $\mathcal E_A=\cspan\{w\1_E:E\in\Sigma,\ E\subseteq A\}$. Since simple functions are dense in $L_2(A_0,a\,\rmd\eta)$, we also have $L_2(A_0,a\,\rmd\eta)=\cspan\{\1_E:E\in\Sigma,\ E\subseteq A_0\}$. Hence $\mathcal E=\cspan\{\1_Ev:E\in\Sigma\}$, so $v$ is cyclic. By Theorem~\ref{thm:representation-james}, the multiplication projections are shrinking in (1) and boundedly complete in (2).
\endgroup
\end{proof}

\begingroup
\begin{remark}\label{rem:associated-zero-set}
The spaces constructed in the proof of Corollary~\ref{cor:associated-spaces} give cyclic dilations of the induced OVM with multiplication projections. The PVM is shrinking in (1) and boundedly complete in (2). For an unconditional Schauder frame $(x_n,f_n)_{n\in\N}$, the coordinate vectors $e_n=\1_{\{n\}}$ form an unconditional Schauder basis of $\mathcal E$ by order continuity. We have $Se_n=x_n$ and $Tx=\sum_nf_n(x)e_n$. Since its coordinate projections are the multiplication projections, the basis is shrinking in (1) and boundedly complete in (2). For unconditional Schauder frames, the equivalences between shrinking and the absence of $\ell_1$, and between bounded completeness and the absence of $c_0$, were proved by Carando, Lassalle and Schmidberg~\cite[Theorems~3.3 and~3.4]{CarandoLassalleSchmidberg2011}. The construction above provides an unconditional associated basis together with the corresponding geometric property.

The summand $L_2(A_0,a\,\rmd\eta)$ ensures that $\mathcal E$ contains a strictly positive function on $\Omega$. Its synthesis map is zero. Since this summand is a Hilbert lattice, adjoining it preserves the absence of $\ell_1$, the KB property, and reflexivity. If $x_t\neq0$ $\eta$-a.e., this summand can be omitted. For a dilation of the induced OVM alone, $\mathcal E_A$ suffices.
\end{remark}
\endgroup

The following two examples use complex Hilbert spaces.

\begin{example}[Hardy--Toeplitz measure]\label{ex:hardy-toeplitz}
Let $m$ be normalized Haar measure on $\mathbb T$, and let $V:\ell_2(\mathbb N_0)\to L_2(m)$ be the isometry $Ve_n(z)=z^n$, with range $H^2$. The canonical phase observable
\[
\Phi(E)=V^*M_{\1_E}V
\]
is a normalized positive operator-valued measure (POVM)~\cite{Pellonpaa2003}. With the multiplication measure on $L_2(m)$, we have a cyclic dilation with vector $\1_{\mathbb T}$. For bounded Borel $b$, the operator $\int b\,\rmd\Phi=V^*M_bV$ is the Toeplitz operator $P_+M_b|_{H^2}$. Since
\[
\inp{\Phi(E)e_n}{e_n}=m(E)\quad(n\geq0)
\]
and $e_n\xrightarrow{w}0$, every value $\Phi(E)$ with $m(E)>0$ is noncompact. By Proposition~\ref{prop:frame-dilation}, $\Phi$ has no continuous Schauder frame representation on $\mathbb T$.
\end{example}

\begin{example}[Cyclic Koopman spectral measures]\label{ex:koopman-cyclic}
Let $Uh=h\circ\tau$ be the Koopman operator of an invertible measure-preserving transformation on a probability space. For its spectral measure $P_U$ and a unit vector $g$, let $H_g=\cspan\{U^kg:k\in\mathbb Z\}$. Define $\rho_g(E)=\inp{P_U(E)g}{g}$ and $P_g(E)=P_U(E)|_{H_g}$. By the spectral theorem, there is a unitary $W_g:L_2(\rho_g)\to H_g$ with $W_g\1_{\mathbb T}=g$ such that
\[
U|_{H_g}=W_gM_zW_g^*\quad\text{and}\quad P_g(E)=W_gM_{\1_E}W_g^*.
\]
Thus $P_g$ is cyclic. It has a continuous Schauder frame representation on $\mathbb T$ exactly when $\rho_g$ is purely atomic. Suppose that a frame representation exists. By Proposition~\ref{prop:frame-dilation}, there is a countable measurable partition $(D_n)$ of $\mathbb T$ such that every $P_g(D_n)$ is compact. Each of these projections has finite-dimensional range, which is unitarily equivalent to $L_2(D_n,\rho_g)$. Let $N=\mathbb T\setminus\{z\in\mathbb T:\rho_g(\{z\})>0\}$, so $\rho_g|_N$ is nonatomic. If $\rho_g(D_n\cap N)>0$, then $L_2(D_n\cap N,\rho_g)$ is infinite-dimensional, a contradiction. Hence $\rho_g(D_n\cap N)=0$ for every $n$, and countable additivity gives $\rho_g(N)=0$. In the atomic case, we have an orthonormal expansion in the normalized characteristic functions. The moments $\int z^k\,\rmd\rho_g=\inp{U^kg}{g}$ are the correlations used in data-driven spectral approximation~\cite{KordaPutinarMezic2020}.
\end{example}

\section{Inductive limits and Barron spaces}\label{sec:inductive}

We now study systems of local cyclic dilations whose analysis, synthesis, and localization operators commute with the connecting maps. These systems recover the norm and topology of the original space. We begin with cyclic subrepresentations of a PVM.

\subsection[Cyclic representations]{Cyclic representations}

\begin{proposition}\label{prop:cyclic-inductive}
Let $P\colon\Sigma\to\cB(X)$ be a PVM. Suppose that every countable subset of $X$ is contained in some $X(u)$. With inherited norms and inclusion maps, these cyclic subspaces form a countably directed family, and
\[
X=\varinjlim_u X(u)\cong\varinjlim_u L_1(m_u)
\]
in the category of Hausdorff locally convex spaces, where $m_u(E)=P(E)u$. Integration defines the second representation, and its connecting maps commute with multiplication by characteristic functions.
\end{proposition}

\begin{proof}
Given $(u_n)$, choose $v$ with $u_n\in X(v)$ for every $n$. Since $X(v)$ is closed and invariant, $X(u_n)\subseteq X(v)$. Thus the family is countably directed and covers $X$. The inclusions are isometries, so the locally convex inductive limit topology is the norm topology of $X$ by~\cite[Theorems~1 and~2]{Pym1968}.

By Proposition~\ref{prop:vector-measure-dilation} with $S=T=I_{X(u)}$, integration and the associated map from $X(u)$ to $L_1(m_u)$ are inverses on simple functions. By density, $I_{m_u}\colon L_1(m_u)\to X(u)$ is an isomorphism. For $X(u)\subseteq X(v)$, the connecting map is $I_{m_v}^{-1}j_{uv}I_{m_u}$, where $j_{uv}$ is inclusion. Since $I_{m_u}(\1_Eh)=P(E)I_{m_u}h$, these maps commute with multiplication by characteristic functions.
\end{proof}

\subsection[Compatible dilations of measure quotients]{Compatible dilations of measure quotients}

Let $\mathcal M(\Omega)$ be the Banach lattice of finite signed countably additive measures on $(\Omega,\Sigma)$, with the total variation norm. For the restriction measure $P(E)\lambda=\lambda|_E$, the cyclic subspace generated by a probability measure $\rho$ is $L_1(\rho)\rho$. Every sequence $(\lambda_n)$ belongs to one such subspace, since we can normalize $\sum_n2^{-n}\abs{\lambda_n}/(1+\norm{\lambda_n}_{\mathrm{TV}})$, unless all the measures vanish. Thus Proposition~\ref{prop:cyclic-inductive} applies. Since $\mathcal M(\Omega)$ is an AL-space, the following lemma is a direct consequence of~\cite[p.~82]{DalesLaustsenOikhbergTroitsky2017}.

\begin{lemma}\label{lem:measure-positive-majorant}
For every bounded real linear operator $R$ on $\mathcal M(\Omega)$, there is a bounded positive operator $V$ such that $ \abs{R\lambda}\leq V\abs{\lambda} $ and $\norm{V}=\norm{R}$.
\end{lemma}

\begingroup
Let $S\colon\mathcal M(\Omega)\to X$ be a bounded surjection with a bounded linear right inverse $T$, and give $X$ the quotient norm
\[
\norm{x}_X=\inf_{S\lambda=x}\norm{\lambda}_{\mathrm{TV}}.
\]
Define $R=TS$ and $\mu(E)x=S((Tx)|_E)$. For $1\leq p<\infty$ and a probability measure $\rho$, define $S_\rho a=S(a\rho)$ and $X_{p,\rho}=S_\rho L_p(\rho)$, with its quotient norm. For probability measures $\rho,\nu$, define $\rho\preceq\nu$ if $\rho\leq c\nu$ for some $c>0$. If $\rho\preceq\nu$, let $i_{\rho\nu}\colon X_{p,\rho}\to X_{p,\nu}$ be inclusion and define $j_{\rho\nu}a=a\,\rmd\rho/\rmd\nu$.

The positive majorant in Lemma~\ref{lem:measure-positive-majorant} allows us to choose compatible cyclic $L_p$-dilations and recover the quotient norm of $X$.

\begin{theorem}\label{thm:measure-inductive-dilation}
Let $S,T$ be the bounded maps above with $ST=I_X$, and let $X$ carry the stated quotient norm. There is a countably directed family $\mathfrak P_R$ of probability measures with the following properties for every $1\leq p<\infty$.
\begin{enumerate}[label=(\roman*),leftmargin=0pt,itemindent=2em]
\item Local dilations. For $\rho\in\mathfrak P_R$, the bounded maps $S_\rho$ and $T_\rho x=\rmd(Tx)/\rmd\rho$, together with multiplication projections on $L_p(\rho)$, give a cyclic dilation of $\mu|_{X_{p,\rho}}$. The dilation space is reflexive for $p>1$.
\item Compatibility. If $\rho,\nu\in\mathfrak P_R$ and $\rho\preceq\nu$, the following diagrams commute for every $E\in\Sigma$.
\[
\begin{CD}
    X_{p,\rho} @>{T_\rho}>> L_p(\rho)\\
    @V{i_{\rho\nu}}VV @VV{j_{\rho\nu}}V\\
    X_{p,\nu} @>{T_\nu}>> L_p(\nu)
\end{CD}
\qquad
\begin{CD}
    L_p(\rho) @>{S_\rho}>> X_{p,\rho}\\
    @V{j_{\rho\nu}}VV @VV{i_{\rho\nu}}V\\
    L_p(\nu) @>{S_\nu}>> X_{p,\nu}
\end{CD}
\qquad
\begin{CD}
    L_p(\rho) @>{M_{\1_E}}>> L_p(\rho)\\
    @V{j_{\rho\nu}}VV @VV{j_{\rho\nu}}V\\
    L_p(\nu) @>{M_{\1_E}}>> L_p(\nu)
\end{CD}
\]
\item Norm recovery.
\begingroup
\[
X=\bigcup_{\rho\in\mathfrak P_R}X_{p,\rho}\quad\text{and}\quad \norm{x}_X=\inf_{\substack{\rho\in\mathfrak P_R\\x\in X_{p,\rho}}}\norm{x}_{p,\rho}.
\]
\endgroup
\item Topology. The locally convex inductive limit $\varinjlim_{\rho\in\mathfrak P_R}X_{p,\rho}$ carries the quotient norm topology of $X$.
\end{enumerate}
\end{theorem}
\endgroup

\begin{proof}
By Lemma~\ref{lem:measure-positive-majorant}, we can take $V\geq0$ with $\abs{R\lambda}\leq V\abs{\lambda}$ and $\norm{V}=\norm{R}$. Define
\[
\mathfrak P_R=\{\rho\in\mathcal M(\Omega)_+:\rho(\Omega)=1,\ V\rho\leq c_\rho\rho\text{ for some }c_\rho<\infty\}.
\]
Given any probability measure $\eta$ and $b>\norm{V}$, let $ v_b=\sum_{n=0}^{\infty}b^{-n}V^n\eta $. Let $k_b=\norm{v_b}_{\mathrm{TV}}$ and $\rho_b=k_b^{-1}v_b$. Since $V$ is positive, the series converges in total variation to a positive measure, with $v_b\geq\eta$ and $Vv_b=b(v_b-\eta)\leq bv_b$. Also, $1\leq\norm{v_b}_{\mathrm{TV}}\leq\sum_{n=0}^{\infty}(\norm{V}/b)^n$. Therefore,
\begin{equation}\label{eq:resolvent-control}
\eta\leq k_b\rho_b,\quad V\rho_b\leq b\rho_b,\quad\text{and}\quad
1\leq k_b\leq\frac{1}{1-\norm{V}b^{-1}}.
\end{equation}
Thus $\rho_b\in\mathfrak P_R$ and $k_b\to1$ as $b\to\infty$. In particular, every probability measure is boundedly dominated by a member of $\mathfrak P_R$. Given probability measures $(\eta_n)$, apply the construction to $\eta=\sum_{n=1}^{\infty}2^{-n}\eta_n$. If $\eta\leq k\rho$ with $\rho\in\mathfrak P_R$, then $\eta_n\leq2^nk\rho$ for every $n$. Thus $\mathfrak P_R$ is countably directed.

Fix $\rho\in\mathfrak P_R$. For bounded $a$, positivity implies $\abs{V(a\rho)}\leq c_\rho\norm{a}_\infty\rho$. Clearly, $V$ preserves $L_1(\rho)\rho$. Its positive density operator $V_\rho$ satisfies
\[
\norm{V_\rho}_{L_1(\rho)\to L_1(\rho)}\leq\norm{V}=\norm{R}\quad\text{and}\quad\norm{V_\rho}_{L_\infty(\rho)\to L_\infty(\rho)}\leq c_\rho.
\]
Since $\abs{R(a\rho)}\leq V(\abs{a}\rho)$, the density $R_\rho a=\rmd R(a\rho)/\rmd\rho$ exists and satisfies $\abs{R_\rho a}\leq V_\rho\abs{a}$. By interpolation and domination, we have
\begin{equation}\label{eq:measure-density-bound}
\norm{R_\rho}_{L_p(\rho)\to L_p(\rho)}\leq\norm{R}^{1/p}c_\rho^{1-1/p}.
\end{equation}
Since $R^2=R$, $SR=S$, and $\ker R=\ker S$, the operator $R_\rho$ is a projection with kernel $\ker S_\rho$. Thus $X_{p,\rho}$ is a Banach quotient. If $S_\rho a=S_\rho a'$, then $R_\rho(a-a')=0$, so $T_\rho(S_\rho a)=R_\rho a$ is well defined. Since $R(a\rho)=T(S_\rho a)$, it agrees with the stated Radon--Nikod\'ym derivative. By taking the infimum over all representatives of $x$, we have
\[
\norm{x}_{p,\rho}\leq\norm{T_\rho x}_{L_p(\rho)}
\leq\norm{R_\rho}\norm{x}_{p,\rho}.
\]
Moreover, $S_\rho T_\rho x=x$. The restriction projections on $\mathcal M(\Omega)$ are strongly countably additive, so $\mu$ is a normalized OVM. For $x\in X_{p,\rho}$, we have $\mu(E)x=S((Tx)|_E)=S_\rho(\1_ET_\rho x)$. Thus $X_{p,\rho}$ is invariant under $\mu(E)$, and we obtain the stated dilation. Multiplication on $L_p(\rho)$ is strongly countably additive and cyclic with vector $\1_\Omega$. The dilation space is reflexive for $1<p<\infty$. If $\rho\leq c\nu$, let $h=\rmd\rho/\rmd\nu$, so $0\leq h\leq c$. For $1\leq p<\infty$, we have
\[
\norm{j_{\rho\nu}a}_{L_p(\nu)}^p=\int_\Omega\abs{a}^ph^{p-1}\,\rmd\rho\leq c^{p-1}\norm{a}_{L_p(\rho)}^p.
\]
The map $j_{\rho\nu}$ is isometric for $p=1$, and in general $\norm{j_{\rho\nu}}\leq c^{1-1/p}$. The same bound holds for $i_{\rho\nu}$ by the quotient norm. Since $(j_{\rho\nu}a)\nu=a\rho$, the synthesis maps agree. For $x\in X_{p,\rho}$, we have
\[
T_\nu i_{\rho\nu}x=\frac{\rmd(Tx)}{\rmd\nu}=\frac{\rmd(Tx)}{\rmd\rho}\frac{\rmd\rho}{\rmd\nu}=(T_\rho x)h=j_{\rho\nu}T_\rho x.
\]
Multiplication by $h$ commutes with multiplication by $\1_E$. This proves commutativity of the three diagrams.

Since $\norm{S_\rho a}_X\leq\norm{a}_{L_1(\rho)}\leq\norm{a}_{L_p(\rho)}$, taking infima gives $\norm{x}_X\leq\norm{x}_{p,\rho}$ for $x\in X_{p,\rho}$. Thus each inclusion $X_{p,\rho}\to X$ is contractive.
\begingroup
For a probability measure $\rho$ and $\lambda\neq0$, let $s=\norm{\lambda}_{\mathrm{TV}}$. If $\abs{\lambda}\leq sc\rho$ for some finite $c$, then $a=\rmd\lambda/\rmd\rho$ satisfies $\norm{a}_{L_1(\rho)}=s$ and $\norm{a}_{L_\infty(\rho)}\leq sc$, so
\[
\norm{a}_{L_p(\rho)}\leq\norm{a}_{L_\infty(\rho)}^{1-1/p}\norm{a}_{L_1(\rho)}^{1/p}\leq sc^{1-1/p}.
\]
Given $S\lambda=x\neq0$, apply~\eqref{eq:resolvent-control} to $\eta=s^{-1}\abs{\lambda}$. Hence $\norm{x}_{p,\rho_b}\leq s k_b^{1-1/p}$. By letting $b\to\infty$ and taking the infimum over $\lambda$, we obtain the norm formula and the union assertion.

It remains to identify the inductive limit topology. Let $q$ be a seminorm continuous on every $X_{p,\rho}$. If $q$ were not continuous on $X$, choose $q(x_n)\geq1$ and $S\lambda_n=x_n$ with $\norm{\lambda_n}_{\mathrm{TV}}\leq2^{-2n}$. The probability measure $\eta=\sum_n2^{-n}\abs{\lambda_n}/\norm{\lambda_n}_{\mathrm{TV}}$ satisfies $\eta\leq k\rho$ for some $\rho\in\mathfrak P_R$ and finite $k$. Since $\abs{\lambda_n}\leq2^nk\norm{\lambda_n}_{\mathrm{TV}}\rho$, the same estimate yields
\[
\norm{x_n}_{p,\rho}\leq\norm{\lambda_n}_{\mathrm{TV}}(2^nk)^{1-1/p}\leq2^{-2n}(2^nk)^{1-1/p}\longrightarrow0,
\]
contradicting continuity of $q|_{X_{p,\rho}}$. Hence the inductive limit topology is the quotient norm topology.
\endgroup
\end{proof}

\subsection[ReLU Barron spaces]{ReLU Barron spaces}\label{sec:barron}

E, Ma and Wu~\cite[Section~2.1]{EMaWu2022} represent their Barron space as a union of reproducing kernel Hilbert spaces indexed by probability measures. Spek, Heeringa, Schwenninger and Brune~\cite[Section~3]{SpekHeeringaSchwenningerBrune2025} represent integral reproducing kernel Banach spaces as unions of $L_p$-quotient spaces and recover their norms by taking infima over probability measures. Bartolucci, De Vito, Rosasco and Vigogna~\cite{BartolucciDeVitoRosascoVigogna2023} decompose functions in weighted ReLU spaces into a part represented by an even measure and an affine function. The canonical measure representation of ReLU Barron functions provides the bounded linear section required by Theorem~\ref{thm:measure-inductive-dilation}. We use it to construct compatible cyclic $L_p$-dilations of a fixed OVM. Modulo affine functions, these dilations describe parameter localization and its exact error. Some of the restricted OVMs admit no continuous Schauder frame representation on the same measurable space.

\begin{example}[ReLU Barron spaces]\label{ex:barron}
Let $\Theta=\mathbb S^d$ where $d\geq1$. Let $\phi_{(w,b)}(x)=(w\cdot x+b)_+$ and define
\[
S\lambda(x)=\int_\Theta\phi_\theta(x)\,\rmd\lambda(\theta)\quad(x\in\R^d).
\]
The Barron space on $\R^d$ is $\mathcal B_d=S\mathcal M(\Theta)$, with norm
\[
\norm{f}_{\mathcal B_d}=\inf_{S\lambda=f}\norm{\lambda}_{\mathrm{TV}}.
\]
We use Euclidean normalization of the parameters. Boundedness of $\phi_\theta(x)$ in $\theta$ shows that $\ker S$ is closed.

\textit{Canonical measure representation.}
\begin{samepage}
Let $\iota(\theta)=-\theta$, $q_\pm=(0,\pm1)$, $u_i^\pm=(\pm e_i,0)$, and $\Theta_0=\Theta\setminus\{q_+,q_-\}$. Denote by $\iota_*$ the pushforward map on $\mathcal M(\Theta)$, defined by $(\iota_*\lambda)(E)=\lambda(\iota^{-1}(E))$ for every Borel set $E\subseteq\Theta$. Since $\phi_{q_+}=1$ and $\phi_{q_-}=0$, the poles contribute only to the affine part. The map $C$ extracts the even measure, $A$ records the affine coefficients, and $J$ represents affine functions by fixed atoms. Define
\begin{align*}
C\lambda&=\frac{1}{2}(\lambda+\iota_*\lambda)|_{\Theta_0},
\quad A\lambda=\left(\frac{1}{2}\int_{\Theta_0}w\,\rmd\lambda,
\frac{1}{2}\int_{\Theta_0}b\,\rmd\lambda+\lambda\{q_+\}\right),
\end{align*}
and $J(v,c)=\sum_{i=1}^dv_i(\delta_{u_i^+}-\delta_{u_i^-})+c\delta_{q_+}$.
\par
\end{samepage}
By the identity $s_+=(\abs{s}+s)/2$, we have
\[
S\lambda=S(C\lambda)+\ell_{A\lambda}\quad\text{and}\quad SJ(v,c)=\ell_{(v,c)},
\]
where $\ell_{(v,c)}(x)=v\cdot x+c$. The identities $C^2=C$, $CJ=AC=0$, and $AJ=I$ show that $R=C+JA$ is a projection with $SR=S$. The Cauchy--Schwarz inequality gives
\[
\norm{J(v,c)}_{\mathrm{TV}}\leq\sqrt{4d+1}\norm{(v,c)}_2\quad\text{and}\quad
\norm{A\lambda}_2\leq\frac{1}{2}\abs{\lambda}(\Theta_0)+\abs{\lambda\{q_+\}},
\]
so $\norm{R}\leq1+\sqrt{4d+1}$.

To identify the kernel, define $\widetilde S\lambda(x,t)=\int_\Theta(w\cdot x+bt)_+\,\rmd\lambda$. E and Wojtowytsch~\cite[Remark~4.6]{EWojtowytsch2022} give
\begin{equation}\label{eq:homogeneous-barron-kernel}
\widetilde S\lambda=0\quad\Longleftrightarrow\quad
\iota_*\lambda=-\lambda\ \text{ and }\ \int_\Theta\theta\,\rmd\lambda(\theta)=0.
\end{equation}
If $S\lambda=0$, homogeneity and continuity give $F:=\widetilde S\lambda=0$ for $t\geq0$, while
\[
F(x,t)-F(-x,-t)=u\cdot x+\beta t,\quad (u,\beta)=\int_\Theta(w,b)\,\rmd\lambda.
\]
Hence $u=0$ and $F=-\beta(-t)_+$. Applying~\eqref{eq:homogeneous-barron-kernel} to $\lambda+\beta\delta_{q_-}$, we obtain $C\lambda=0$, and the synthesis identity gives $A\lambda=0$. Together with $SR=S$, this proves $\ker R=\ker S$. Thus each $f\in\mathcal B_d$ has a unique representation
\[
f(x)=\frac{1}{2}\int_{\Theta_0}\abs{w\cdot x+b}\,\rmd\gamma_f(w,b)+v_f\cdot x+c_f,
\]
where $\gamma_f$ is a finite signed measure on $\Theta_0$ satisfying $\iota_*\gamma_f=\gamma_f$, and $(\gamma_f,(v_f,c_f))=(C\lambda,A\lambda)$ for any $S\lambda=f$. Thus $Tf=\gamma_f+J(v_f,c_f)$ defines a linear right inverse of $S$, with $T(S\lambda)=R\lambda$ and $\norm{T}\leq1+\sqrt{4d+1}$.

\Needspace{6\baselineskip}
\textit{Compatible cyclic dilations.} Let $\mathcal A=\{u_i^+,u_i^-:1\leq i\leq d\}\cup\{q_+,q_-\}$. Define
\[
\begin{aligned}
    \mathfrak P=\{\rho\in\mathcal M(\Theta)_+:{}&\rho(\Theta)=1,\ \iota_*\rho\leq k_\rho\rho\text{ for some }k_\rho<\infty,\\
    &\rho\{u\}>0 \text{ for all } u\in\mathcal A\}.
\end{aligned}
\]
For $\rho\in\mathfrak P$, let $r_\rho=\rmd(\iota_*\rho)/\rmd\rho$. The reflection condition and the positive atomic masses give $R(a\rho)\ll\rho$, with density
\begin{align*}
R_\rho a
=\frac{1}{2}\1_{\Theta_0}(a+r_\rho a\circ\iota)+\sum_{i=1}^dv_i(a)\left(\frac{\1_{\{u_i^+\}}}{\rho\{u_i^+\}}-\frac{\1_{\{u_i^-\}}}{\rho\{u_i^-\}}\right)
  +c(a)\frac{\1_{\{q_+\}}}{\rho\{q_+\}},
\end{align*}
where $(v(a),c(a))=A(a\rho)$. Let $m_\rho=\min_{u\in\mathcal A}\rho\{u\}>0$. The density identity and the displayed formula show that
\[
\norm{R_\rho}_{L_1\to L_1}\leq1+\sqrt{4d+1},
\]
and
\[
\norm{R_\rho}_{L_\infty\to L_\infty}\leq\frac{1+k_\rho}{2}+\frac{1}{2m_\rho}    \leq\frac{1+k_\rho}{2}+\frac{\sqrt{4d+1}}{m_\rho}=:D_\rho.
\]
Here $k_\rho\geq1$, and the correction terms are supported on distinct atoms. Interpolation yields
\[
\norm{R_\rho}_{L_p\to L_p}\leq(1+\sqrt{4d+1})^{1/p}D_\rho^{1-1/p}\quad(1\leq p<\infty).
\]
In the notation of Theorem~\ref{thm:measure-inductive-dilation}, denote $X_{p,\rho}$ by $\mathcal B_{p,\rho}$. The identities $R^2=R$ and $\ker R=\ker S$ imply $R_\rho^2=R_\rho$ and $\ker R_\rho=\ker S_\rho$. The projection and compatibility arguments in the proof of Theorem~\ref{thm:measure-inductive-dilation} therefore give
\[
f\in\mathcal B_{p,\rho}\quad\Longleftrightarrow\quad
Tf\ll\rho\ \text{ and }\ a_f:=\frac{\rmd(Tf)}{\rmd\rho}\in L_p(\rho).
\]
The density $a_f$ is the unique coefficient fixed by $R_\rho$ with $S_\rho a_f=f$, and $ \norm{f}_{p,\rho}\leq\norm{a_f}_p\leq\norm{R_\rho}\norm{f}_{p,\rho}$. With $T_\rho f=a_f$, the compatible cyclic dilations of $\mu(E)f=S((Tf)|_E)$ are
\[
\mu(E)|_{\mathcal B_{p,\rho}}=S_\rho M_{\1_E}T_\rho=:\mu_\rho(E)\quad\text{and}\quad S_\rho T_\rho=I.
\]
The dilation space $L_p(\rho)$ is reflexive for $1<p<\infty$. For $\rho\preceq\nu$, the analysis, synthesis, and multiplication operators commute with the connecting maps in Theorem~\ref{thm:measure-inductive-dilation}.

Let $\tau$ be uniform on $\mathcal A$. For every probability measure $\eta$ and $0<\varepsilon<1$, let $ \rho_\varepsilon=(1-\varepsilon)\eta+\frac{\varepsilon}{2}(\iota_*\eta+\tau) $. Reflection exchanges $\eta$ and $\iota_*\eta$, and $\tau$ is symmetric with positive mass on $\mathcal A$, so $\rho_\varepsilon\in\mathfrak P$ and $\eta\leq(1-\varepsilon)^{-1}\rho_\varepsilon$. For $(\rho_n)\subseteq\mathfrak P$, applying this construction to $\eta=\sum_n2^{-n}\rho_n$ with $\varepsilon=1/2$ gives a common upper bound. Thus $\mathfrak P$ is countably directed, and $(1-\varepsilon)^{-1}\to1$ allows us to apply the norm and topology arguments in the proof of Theorem~\ref{thm:measure-inductive-dilation}:
\[
\mathcal B_d=\varinjlim_{\rho\in\mathfrak P}\mathcal B_{p,\rho}\quad\text{and}\quad
\norm{f}_{\mathcal B_d}=\inf_{\substack{\rho\in\mathfrak P\\f\in\mathcal B_{p,\rho}}}\norm{f}_{p,\rho},
\]
with the Barron norm topology. These quotient spaces have continuous point evaluations and are reproducing kernel Hilbert spaces when $p=2$.

\textit{Quotient localization.} The quotient description gives the exact error when the nonlinear part is restricted to a symmetric parameter set and the affine part is retained.
\begingroup
Let $\Sigma_{\mathrm{sym}}$ be the sigma-algebra of symmetric Borel subsets of $\Theta$, and let $M_E\lambda=\lambda|_E$. For $E\in\Sigma_{\mathrm{sym}}$, both $\operatorname{ran}C$ and $\operatorname{ran}J$ are invariant under $M_E$, so $RM_ER=M_ER$. Thus $\mu(E)\mu(F)=SM_ERM_FT=\mu(E\cap F)$ for $E,F\in\Sigma_{\mathrm{sym}}$, and $\mu|_{\Sigma_{\mathrm{sym}}}$ is a normalized PVM. Factoring out affine functions gives the following isometric representation.
\par\endgroup
Let $\operatorname{Aff}_d$ denote the space of affine functions on $\R^d$, and let $q\colon\mathcal B_d\to\mathcal B_d/\operatorname{Aff}_d$ be the quotient map. Contractivity of $C$ and the representative $S\gamma_f$ give
\begin{equation}\label{eq:barron-even-isometry}
\norm{qf}_{\mathcal B_d/\operatorname{Aff}_d}
=\norm{\gamma_f}_{\mathrm{TV}}
=\norm{S\gamma_f}_{\mathcal B_d}.
\end{equation}
Hence $qf\mapsto\gamma_f$ is an onto isometry to $\mathcal M_{\mathrm{even}}(\Theta_0)$, the space of finite even signed measures on $\Theta_0$. For a symmetric Borel set $E\subseteq\Theta_0$, define $\widehat P(E)qf=qS(\gamma_f|_E)$. Under this isometry, $\widehat P$ is the restriction PVM on the sigma-algebra of symmetric Borel subsets of $\Theta_0$.
\begingroup
Since $C(\lambda|_E)=(C\lambda)|_E$, the following diagrams commute for every $\rho\in\mathfrak P$.
\begin{equation}\label{eq:barron-quotient-dilation}
\begin{CD}
    \mathcal B_d @>{\mu(E)}>> \mathcal B_d\\
    @V{q}VV @VV{q}V\\
    \mathcal B_d/\operatorname{Aff}_d @>{\widehat P(E)}>> \mathcal B_d/\operatorname{Aff}_d
\end{CD}
\qquad
\begin{CD}
    L_p(\rho) @>{M_{\1_E}}>> L_p(\rho)\\
    @V{qS_\rho}VV @VV{qS_\rho}V\\
    \mathcal B_d/\operatorname{Aff}_d @>{\widehat P(E)}>> \mathcal B_d/\operatorname{Aff}_d
\end{CD}
\end{equation}
\endgroup
Thus the multiplication projections in the cyclic dilations of $\mu$ induce $\widehat P$ on the corresponding quotient spaces $\mathcal B_{p,\rho}/\operatorname{Aff}_d$.
\begingroup
Let $H_d=\Theta_0/\{\theta\sim-\theta\}$, with its quotient Borel sigma-algebra, and let $\varpi\colon\Theta_0\to H_d$ be the quotient map. For $1\leq p<\infty$ and $\rho\in\mathfrak P$, define the finite measure $\nu_\rho=\varpi_*(\rho|_{\Theta_0})$. The positive masses on $\mathcal A$ imply $\operatorname{Aff}_d\subseteq\mathcal B_{p,\rho}$. Let $\mathcal Q_{p,\rho}=\mathcal B_{p,\rho}/\operatorname{Aff}_d$, and let $q_\rho$ be the quotient map. Then
\[
W_\rho\colon\mathcal Q_{p,\rho}\longrightarrow L_p(\nu_\rho) \qquad\text{where}\qquad
W_\rho q_\rho f=\frac{\rmd(\varpi_*\gamma_f)}{\rmd\nu_\rho}=:g_f
\]
is an onto isometry. Since pushforward by $\varpi$ is injective on even measures, $W_\rho$ is injective. If $f=S_\rho a$, then $\varpi_*\gamma_f=\varpi_*(a\rho|_{\Theta_0})$, and Jensen's inequality gives $\norm{g_f}_{L_p(\nu_\rho)}\leq\norm{a}_{L_p(\rho)}$. Taking infima over representations and affine perturbations gives $\norm{g_f}_{L_p(\nu_\rho)}\leq\norm{q_\rho f}_{\mathcal Q_{p,\rho}}$. Conversely, let $g\in L_p(\nu_\rho)$ and define $a=\1_{\Theta_0}g\circ\varpi$, with $a=0$ at the poles. Then $W_\rho q_\rho S_\rho a=g$ and $\norm{a}_{L_p(\rho)}=\norm{g}_{L_p(\nu_\rho)}$. Under this isometry, the induced projection $\widehat P_\rho(E)$ on $\mathcal Q_{p,\rho}$ is multiplication by $\1_{\varpi(E)}$ for every symmetric Borel set $E\subseteq\Theta_0$. In particular,
\[
\norm{q_\rho f-\widehat P_\rho(E)q_\rho f}_{\mathcal Q_{p,\rho}}
=\left(\int_{H_d\setminus\varpi(E)}\abs{g_f}^p\,\rmd\nu_\rho\right)^{1/p}.
\]
\par\endgroup

Every countable family of finite even measures is dominated by one even probability measure $\nu$. The even part of $L_1(\nu)\nu$ is cyclic with generator $\nu$, so Proposition~\ref{prop:cyclic-inductive} also gives the quotient as an inductive limit of these cyclic subspaces.

Since $(w,b)$ and $(-w,-b)$ define the same hyperplane, $\widehat P(E)$ retains the nonlinear components on the unoriented hyperplanes represented by $E$. Retaining the affine part of $f$, let $f_E=S(\gamma_f|_E)+\ell_{(v_f,c_f)}$. By~\eqref{eq:barron-quotient-dilation}, we have $qf_E=\widehat P(E)qf=q\mu(E)f$. Equation~\eqref{eq:barron-even-isometry} applied to $\gamma_f|_{\Theta_0\setminus E}$ gives
\begin{equation}\label{eq:barron-localization-error}
\begin{aligned}
    \norm{qf-\widehat P(E)qf}_{\mathcal B_d/\operatorname{Aff}_d}=\norm{f-f_E}_{\mathcal B_d}=\abs{\gamma_f}(\Theta_0\setminus E).
\end{aligned}
\end{equation}
The Cauchy--Schwarz inequality also gives
\[
\abs{f(x)-f_E(x)}\leq\frac{1}{2}\sqrt{1+\norm{x}_2^2}\,\abs{\gamma_f}(\Theta_0\setminus E).
\]
For a finite network, this operation removes the components on the discarded hyperplanes after combining their coefficients and retains the affine term. If $\gamma_f\neq0$, the probability measure $\nu_f=\abs{\gamma_f}/\norm{\gamma_f}_{\mathrm{TV}}$ is determined by $f$, and its mass on $\Theta_0\setminus E$ is the error in~\eqref{eq:barron-localization-error} divided by $\norm{\gamma_f}_{\mathrm{TV}}$.

For the homogeneous model on $\R^d$, equivalently on its input sphere, take $\Theta=\mathbb S^{d-1}$ and $\phi_w(x)=(w\cdot x)_+$. Define
\[
C\lambda=\frac{1}{2}(\lambda+\iota_*\lambda),\quad
A\lambda=\frac{1}{2}\int_\Theta w\,\rmd\lambda(w),\quad\text{and}\quad
Jv=\sum_{i=1}^dv_i(\delta_{e_i}-\delta_{-e_i}).
\]
By~\eqref{eq:homogeneous-barron-kernel}, the same formulas apply with the pole terms omitted and affine functions replaced by linear functions.

\textit{Failure of frame representation.}
\begingroup
The specified OVM need not have a continuous Schauder frame representation on the same measurable space. In the biased model, let $\sigma$ be uniform on $\mathbb S^d$ and $\rho=(\sigma+\tau)/2$. Since $\rho$ is symmetric, the map $U_\rho g=S_\rho(\1_{\Theta_0}g\circ\varpi)$, with the coefficient zero at the poles, is an isometric right inverse of $W_\rho q_\rho$. For any Borel set $E\subseteq\Theta$, define
\[
k_E(\varpi\theta)=\frac{\1_E(\theta)+\1_E(-\theta)}{2}\quad(\theta\in\Theta_0).
\]
The dilation formula gives
\[
W_\rho q_\rho\mu_\rho(E)U_\rho=M_{k_E}.
\]
If $\sigma(E)>0$, then $k_E\geq\frac{1}{2}$ on a set of positive $\nu_\rho$-measure disjoint from its atoms. Thus $M_{k_E}$ is noncompact, and the factorization shows that $\mu_\rho(E)$ is noncompact. Since $\mathcal B_{p,\rho}$ is separable, Proposition~\ref{prop:frame-dilation} shows that a continuous Schauder frame representation would give a countable cover of $\Theta$ by sets with compact operator values. Every set in this cover would be $\sigma$-null, a contradiction.
\par\endgroup
For $p=2$, the space $\mathcal B_{2,\rho}$ is a separable Hilbert space and therefore has an orthonormal basis, although the OVM $\mu_\rho$ has no continuous Schauder frame representation on $\Theta$.
\end{example}

\section{Examples}\label{sec:examples}

The following examples separate shrinking from bounded completeness and show that the representation results also apply to spaces not isomorphic to Banach lattices.

\begin{example}[Boundedly complete but not shrinking]\label{ex:l1-james}
On $X=L_1[0,1]$, the multiplication measure $P(E)h=\1_Eh$ is cyclic with vector $\1$. For $m(E)=\1_E$, we have $L_1(m)=L_1[0,1]$ isometrically and $I_m=I_X$. Since disjointly supported vectors satisfy $ \norm{\sum_{n=1}^Nz_n}_1=\sum_{n=1}^N\norm{z_n}_1$, the measure is boundedly complete. It is not shrinking, since $(0,1/n)\downarrow\varnothing$ while $\norm{P((0,1/n))^*\1}_\infty=1$. The space admits neither a continuous nor an unconditional discrete Schauder frame: it contains no copy of $c_0$ and fails the Radon--Nikod\'ym property, since $m$ has Lebesgue variation and no Bochner density~\cite[Chapter~III]{DiestelUhl1977}, whereas Propositions~\ref{prop:frame-dilation} and~\ref{prop:variation-rnp} would imply this property.
\end{example}

\begin{example}[Shrinking but not boundedly complete]\label{ex:orlicz-james}
Let $\Phi(u)=e^{u^2}-1$ for $u\geq0$, and consider the Orlicz heart
\[
G=(\exp L_2)_0=\overline{L_\infty[0,1]}^{\,L_\Phi},
\]
with the Luxemburg norm $\norm{h}_\Phi=\inf\{\lambda>0:\int_0^1\Phi(\abs{h(t)}/\lambda)\,\rmd t\leq1\}$. The space $G$ is separable and has order-continuous norm. Hence $P(E)h=\1_Eh$ is strongly countably additive and cyclic with vector $\1$. By Orlicz duality~\cite{LindenstraussTzafriri1979}, the dual is $L_\Psi$ with its Orlicz norm, where $\Psi$ is complementary to $\Phi$. Since $\Psi(v)\asymp v\sqrt{\log v}$ as $v\to\infty$, it satisfies the $\Delta_2$-condition at infinity, and $L_\Psi$ has order-continuous norm. Since $P(E)^*g=\1_Eg$, the adjoint measure is strongly countably additive. Thus $P$ is shrinking.

Choose pairwise disjoint $A_n\subseteq[0,1]$ with $\abs{A_n}=e^{-n^2}$, and let $z_n=n\1_{A_n}$. Then $z_n\in P(A_n)G$, and
\begin{align*}
\int_0^1\Phi\left(\frac{1}{2}\sum_{n=1}^Nz_n(t)\right)\,\rmd t
=\sum_{n=1}^Ne^{-n^2}\bigl(e^{n^2/4}-1\bigr)<1.
\end{align*}
Hence $\norm{\sum_{n=1}^Nz_n}_\Phi\leq2$, whereas $\norm{z_n}_\Phi=n/\sqrt{\log(1+e^{n^2})}\to1$. Thus $P$ is not boundedly complete.

The space $G$ admits neither an unconditional discrete Schauder frame nor a continuous Schauder frame. To prove this, first consider the Haar functions
\[
h_k=\1_{[0,2^{-k-1})}-\1_{[2^{-k-1},2^{-k})}\quad(k\geq0).
\]
The alternating sums are bounded in absolute value by $2$, whereas the sum with all positive coefficients equals $N$ on $[0,2^{-N})$. Therefore
\begin{align*}
\norm{\sum_{k=0}^{N-1}(-1)^kh_k}_\Phi
\leq2\norm{\1}_\Phi\quad\text{and}\quad
\norm{\sum_{k=0}^{N-1}h_k}_\Phi
\geq\frac{N}{\sqrt{\log(1+2^N)}}\longrightarrow\infty.
\end{align*}
Thus the Haar basis of $G$ is conditional. By the Haar reproducibility theorem~\cite[Theorem~2.c.8]{LindenstraussTzafriri1979}, a separable symmetric function space which embeds into a space with an unconditional basis has an unconditional Haar basis. Hence $G$ admits no such embedding, and therefore no unconditional Schauder frame by~\cite[Proposition~1.6]{HanLarsonLiuLiu2014}. By Corollary~\ref{cor:frame-discretization}, a continuous Schauder frame would imply an unconditional Schauder frame. Thus $G$ admits no continuous Schauder frame.
\end{example}

\begin{example}[Representations on spaces not isomorphic to Banach lattices]\label{ex:vector-sequences}
Let $Y$ be the separable uniformly convex HI space of Ferenczi~\cite{Ferenczi1997}. Consider the coordinate measure $(P(A)y)_n=\1_A(n)y_n$ on $c_0(Y)$, $\ell_1(Y)$, and $\ell_2(Y)$. For $F_m=\{1,\ldots,m\}$, the range of $P(F_m)$ is reflexive, so $P(F_m)^{**}X^{**}\subseteq J_XX$. Thus $P$ satisfies the localization property~\eqref{eq:localization}. These coordinate representations are not of finite multiplicity. Indeed, if a nonzero spectral subspace $P(A)X$ were generated by $r$ vectors $x_1,\ldots,x_r$, then any coordinate $n\in A$ would give $Y=\spn\{x_1(n),\ldots,x_r(n)\}$, contrary to the infinite dimensionality of $Y$. On $c_0(Y)$, the adjoint acts on $\ell_1(Y^*)$, so $P$ is shrinking, while single-coordinate unit vectors show failure of bounded completeness. On $\ell_1(Y)$, additivity of the norm on disjoint supports implies bounded completeness, whereas a constant unit functional in $\ell_\infty(Y^*)$ shows failure of shrinking. On $\ell_2(Y)$, reflexivity implies both properties. A nonatomic example is multiplication on $L_2([0,1];Y)$. By Proposition~\ref{prop:bochner-normal}, this representation also satisfies~\eqref{eq:localization}, and its measure is both shrinking and boundedly complete.

Each of these spaces contains a complemented copy of $Y$. Since $Y$ contains no unconditional basic sequence, it is not isomorphic to a complemented subspace of a Banach lattice~\cite[Section~1]{DeHeviaTradacete2025}. Thus none of these spaces is isomorphic to a Banach lattice. An unconditional Schauder frame would complement the space in a space with an unconditional basis~\cite[Proposition~1.6]{HanLarsonLiuLiu2014}, and hence complement $Y$ in a Banach lattice after an equivalent renorming. Therefore these spaces admit neither unconditional nor continuous Schauder frames by Corollary~\ref{cor:frame-discretization}.
\end{example}

\begingroup
The properties established above are summarized below, with $Y$ as in Example~\ref{ex:vector-sequences}.
\begin{center}
\small
\begin{tabular}{lcccc}
\hline
$X$ & \shortstack{Shrinking\\$P$} & \shortstack{Boundedly\\complete $P$} & \shortstack{Finite\\multiplicity} & \shortstack{Continuous\\frame on $X$}\\
\hline
$L_1[0,1]$ & No & Yes & Yes & No\\
$(\exp L_2)_0$ & Yes & No & Yes & No\\
$c_0(Y)$ & Yes & No & No & No\\
$\ell_1(Y)$ & No & Yes & No & No\\
$\ell_2(Y)$ & Yes & Yes & No & No\\
\hline
\end{tabular}
\end{center}
\par\endgroup

\section*{Acknowledgments}

Daniel Freeman was supported by the U.S. National Science Foundation (Awards DMS~2154931 and DMS~2556024).
Qianfeng Hu acknowledges support from the National Natural Science Foundation of China (Grant Nos.~12301162 and~12671162), the Hebei Natural Science Foundation (Grant No.~A2024202039), and Young Scientific and Technological Talents (Level Three) in Tianjin (Grant No.~QN20250339).
Rui Liu and Jie Shen were partially supported by the National Natural Science Foundation of China (Grant Nos.~12071230 and~12471131).

\section*{Declaration of AI use}

\begin{enumerate}[itemsep=6pt]
\item \textbf{Scope of use.} The authors used GPT 6 (OpenAI) which is the chat model, and GPT 6 Astra (OpenAI) during the research and preparation of this manuscript. The tools assisted with mathematical discussions, drafting and revision, preparing commutative diagrams and summary tables, checking statements and proofs, and identifying potentially relevant prior work and questions of novelty.

\item \textbf{AI contributions to mathematical ideas.} The AI tools suggested applying cyclicity, weak compactness techniques, and Banach lattice methods to specific problems concerning PVMs in this paper. The initial suggestions for Remark~\ref{rem:lattice-james}, Lemma~\ref{lem:multiplier-hull}, Theorem~\ref{prop:multiplier-orbit-criterion}, and Examples~\ref{ex:no-reflexive-dilation}, \ref{ex:orbit-variation-no-frame}, \ref{ex:c0-variation-no-frame}, \ref{ex:hardy-toeplitz}, \ref{ex:koopman-cyclic}, and~\ref{ex:vector-sequences} came from these tools. Further discussions with the AI tools contributed to the development or strengthening of Lemma~\ref{lem:normal-restriction}, the three-space property of the localization condition in Proposition~\ref{prop:normal-coverage}, Remark~\ref{rem:normal-quotients}, the removal of the separability assumption in Proposition~\ref{prop:bochner-normal}, the simplification of the proof of Theorem~\ref{thm:dilation-geometry} using the DFJP construction, the applications to Barron spaces in Section~\ref{sec:barron}, and Example~\ref{ex:orlicz-james}.

\item \textbf{Ideas originating from the authors.} The authors' principal research ideas include the localization property and its equivalent characterizations; the connection between James-type theorems for operator-valued measures and dilation theory; the use of K\"othe function spaces in this setting; the characterization of reflexive dilations under local variation conditions (Theorem~\ref{thm:variation-reflexive-dilation}); the relationship between continuous Schauder frames and operator-valued measures; and the construction of compatible cyclic $L_p$-dilations, including the inductive-limit framework of Theorem~\ref{thm:measure-inductive-dilation} and its applications to Barron spaces.

\item \textbf{Verification and responsibility.} The authors have carefully reviewed and verified all AI-assisted mathematical statements, arguments, constructions, and examples incorporated into this paper, revising or rewriting them where necessary. AI-assisted comparisons with previous work were treated as preliminary suggestions rather than independent verification of novelty. The authors made the final decisions on mathematical correctness, attribution, and presentation, and take full responsibility for the arguments, references, claims of originality, and final manuscript.
\end{enumerate}

\section*{Lean 4 formal verification}

The mathematical results specified below were formalized and their proofs were checked by the Lean~4 kernel. The source code and verification materials are available in the accompanying GitHub repository \href{https://github.com/JoesephShen/The-Geometry-of-Banach-space-operator-valued-measures-with-cyclic-dilations-Lean-4}{ https://github.com/JoesephShen/The-Geometry-of-Banach-space-operator-valued-measures-with-cyclic-dilations-Lean-4}.

\begin{enumerate}[itemsep=6pt,before={\emergencystretch=2em}]
\item \textbf{Scope.} The formalization covers the three main theorems in the introduction and all numbered mathematical results in the body of the paper, including theorems, lemmas, propositions, corollaries, remarks, and examples.

\item \textbf{Axiom dependencies.} The verification uses the standard Lean axioms \texttt{propext}, \texttt{Classical.choice}, and \texttt{Quot.sound}, as well as 70 textbook and literature results declared explicitly as axioms in the Lean source. The verified proofs are checked relative to these assumptions. The repository lists the sources of the external results and records the transitive axiom dependencies of 620 audited declarations, none of which depends on \texttt{sorryAx}.

\item \textbf{AI assistance.} GPT 6.2 Sol (OpenAI) assisted with generating and reviewing the Lean~4 code.
\end{enumerate}

\bibliographystyle{plain}
\bibliography{cite}

\end{document}